\documentclass[]{elsarticle}
\usepackage{subfigure}

\usepackage{amssymb,amsmath,newlfont,enumerate}

\let\eps\varepsilon

\let\ds\displaystyle
\def\qquad{{\quad\quad}}

\usepackage{amsfonts}
\usepackage{amsbsy}

\usepackage[]{graphicx}
\usepackage[hang,center]{caption}
\usepackage{verbatim} 
\usepackage{float}
\usepackage{color}
\usepackage{multirow}

\begin{document}

%% \title[AP-scheme for the anisotropic parabolic
%% equations]{Asymptotic Preserving scheme for strongly anisotropic
%% parabolic equations for arbitrary anisotropy direction}

\title{Asymptotic Preserving scheme for strongly anisotropic
  parabolic equations for arbitrary anisotropy direction}

\author[J. Narski]{Jacek Narski} 
\address[J. Narski]{Universit\'e de Toulouse,
  Institut de Math\'ematiques de Toulouse, 118 route de Narbonne,
  F-31062 Toulouse, France}

\author[M. Ottaviani]{Maurizio Ottaviani}
\address[M. Ottaviani]{CEA, IRFM, F-13108 Saint-Paul-lez-Durance, France}

\date{\today}

\begin{abstract}
  This paper deals with the numerical study of a strongly anisotropic
  heat equation. The use of standard schemes in this situation leads
  to poor results, due to the high anisotropy.  Furthermore, the
  recently proposed Asymptotic-Preserving method \cite{LNN} allows one
  to perform simulations regardless of the anisotropy strength but its
  application is limited to the case, where the anisotropy direction
  is given by a field with all field lines open. In this paper we
  introduce a new Asymptotic-Preserving method, which overcomes those
  limitations without any loss of precision or increase in the
  computational costs.  The convergence of the method is shown to be
  independent of the anisotropy parameter $0 < \eps <1$, and this for
  fixed coarse Cartesian grids and for variable anisotropy
  directions. The context of this work is magnetically confined fusion
  plasmas.
\end{abstract}


\maketitle
%%%%%%%%%%%%%%%%%%%%%%%%%%%%%%%%%%%%%%%%%%%%%%%%%%%%%%%%%%%
%%%%%%%%%%%%%%%%%%%%%%%%%%%%%%%%%%%%%%%%%%%%%%
\section{Introduction} \label{SEC1}
%%%%%%%%%%%%%%%%%%%%%%%%%%%%%%%%%%%%%%%%%%%%%%%

This work deals with the efficient numerical treatment of heat
transport in a strongly anisotropic medium. We address in particular
models of magnetised plasma with magnetic field perturbations such as
those produced by tearing modes and magnetic islands.

In classical transport theory of strongly magnetised plasmas, the
ratio of the parallel ($\chi_\parallel$) to the perpendicular
($\chi_\perp$)heat conductivity of a given species (electrons or ions)
scales like $(\Omega_c \tau_c)^2$ where $\Omega_c$ is the cyclotron
frequency (the rotation frequency around the field lines) and $\tau_c$
the collision frequency. This product is several orders of magnitude
(typically 10 to 12).

Magnetic islands are non-ideal deformations of the primary magnetic
field. In plasma confinement devices, they have a small magnetic
component pointing outwards. However, due to the strong parallel
conductivity, even a tiny outward components leads to a substantial
heat loss in the island regions. Thus, magnetic islands are unwanted
effects in actual applications.

Theories of the formation of magnetic islands rely on various
ingredients. In the regime where tearing modes (TM) are linearly
unstable, magnetic islands are the result of TM evolution and
saturation. When however TM are stable, magnetic islands can still
occur through a mechanism of self-sustainment. In this regime, a key
element of the island dynamics is the competition between the parallel
and the perpendicular heat fluxes, depending in particular on the
ratio $\chi_\parallel/\chi_\perp$ (denoted in the sequel as
$1/\varepsilon $), which may ultimately determine whether the island
grows or is suppressed.

The heat equation studied in this paper can be written as 
\begin{gather}
  \partial_t u - {1 \over \eps} \nabla_{\parallel} \cdot (A_{\parallel}
  \nabla_{\parallel} u) -\nabla_{\perp} \cdot (A_{\perp} \nabla_{\perp}
  u)=0\,,
  \label{eq:Jqeb}
\end{gather}
where $\nabla_\parallel $ and $\nabla_\perp $ denote the gradient in
the direction parallel (respectively perpendicular) to the magnetic
field. This problem can become (depending on boundary conditions)
ill-posed in the limit of $\varepsilon \rightarrow 0$.

Conventional numerical methods usually fail to give accurate results
with reasonable computational resources for realistic physical
parameters. Indeed, when $1/\eps$ is of the order of $10^{10}$ the
problem becomes extremely anisotropic and standard discretizations
lead to very badly conditioned linear systems with condition number
proportional to $1/(\varepsilon h^2)$ ($h$ being the spatial
discretization step). It is therefore important to develop a numerical
scheme that can address the problem and give accurate results
independently of $\varepsilon $.

The original motivation of this work comes from the fusion plasma
physics, but similar anisotropic problems are encountered in many
other fields of application. One can mention for example image
processing \cite{Perona,Weickert}, transport modeling in fractured
geological structures \cite{Berkowitz} or semiconductor modeling
\cite{semicond}.

Numerical resolution of strongly anisotropic problems has been
addressed by many authors. For example adapted coordinates are often
used in the context of plasma simulation \cite{Beer,Boozer,Hamada}.
This approach can be however difficult to implement, especially when
the magnetic field is variable in time. This is why it is preferable
to choose a method which does not require mesh or coordinate
adaptation, like in \cite{Ottaviani}. Another approach relies on
numerical schemes specially developed for the anisotropic
context. Finite difference schemes were investigated in
\cite{Gunter05,Sharma07,sharma2011fast}. High order finite element
method was proposed in \cite{G}. Multigrid methods \cite{Gee} can
sometimes be beneficial. These methods are usually efficient for a
selected range of $\eps $ but do not behave well in the limit $\eps
\rightarrow 0$.

A different way to overcome this difficulty (adopted in this paper) is
to apply the so called Asymptotic Preserving scheme introduced first
in \cite{ShiJin} to deal with singularly perturbed kinetic models. The
idea is to reformulate the initial problem into an equivalent form,
which remains well-posed, even if the anisotropy strength is
infinite. The reformulation that is studied in this paper was first
applied to the anisotropic stationary diffusion equation in
\cite{DLNN} and then to the nonlinear anisotropic heat equation in
\cite{MN,LNN}. This method is based on introduction of an auxiliary
variable, which serves to eliminate from the equation the dominant
part, {\it i.e.} the one multiplied by $1/\varepsilon $. The choice of
the auxiliary variable presented in those papers allowed to solve the
problem regardless of the anisotropy strength but imposed serious
limitations on the magnetic field. In particular, the case of magnetic
islands cannot be treated by those schemes. In this paper we propose a
new method which overcomes this limitation.

The plan of the article is as follows. In Sec.~\ref{SEC2} the
mathematical problem is presented. Sec.~\ref{SEC3} is devoted to the
description of the numerical method. The numerical tests with known
analytic solutions and the application to the problem of transport in
a magnetic island are presented in Sec~\ref{sec:numres}.

%%%%%%%%%%%%%%%%%%%%%%%%%%%%%%%%%%%%%%%%%%%%%%%
\section{Description of the mathematical problem} \label{SEC2}
%%%%%%%%%%%%%%%%%%%%%%%%%%%%%%%%%%%%%%%%%%%%%%%

We are interested in a resolution of an anisotropic, two or three
dimensional heat problem defined on a domain $\Omega $.  Let the
anisotropy direction be given by a smooth and normalized vector field
$b$, $|b|=1$ and let the computational domain $\Omega $ be a bounded
and sufficiently smooth two or three dimensional subset of $\mathbb
R^d$ with $d=2,3$. The domain $\Omega $ is equipped with a boundary
$\Gamma $, which is decomposed accordingly to the boundary conditions
into two parts: : $\Gamma _D$ and $\Gamma _N= \partial\Omega \setminus
\Gamma _D$ with the Dirichlet and Neumann boundary condition imposed
respectively.

It is convenient to decompose vectors $v \in \mathbb R^d$, gradients
$\nabla \phi$, with $\phi(x)$ a scalar function, and divergences
$\nabla \cdot v$, with $v(x)$ a vector field , into a part parallel to
the anisotropy direction $b$ and a part perpendicular to it. These
parts are defined as follows:
\[
\begin{array}{llll}
  \ds v_{\parallel}:= (v \cdot b) b \,, & \ds v_{\perp}:= (Id- b \otimes b) v\,, &\textrm{such that}&\ds
  v=v_{\parallel}+v_{\perp}\,,\\[3mm]
  \ds \nabla_{\parallel} \phi:= (b \cdot \nabla \phi) b \,, & \ds
  \nabla_{\perp} \phi:= (Id- b \otimes b) \nabla \phi\,, &\textrm{such that}&\ds
  \nabla \phi=\nabla_{\parallel}\phi+\nabla_{\perp}\phi\,,\\[3mm]
  \ds \nabla_{\parallel} \cdot v:= \nabla \cdot v_{\parallel}  \,, & \ds
  \nabla_{\perp} \cdot v:= \nabla \cdot v_{\perp}\,, &\textrm{such that}&\ds
  \nabla \cdot v=\nabla_{\parallel}\cdot v+\nabla_{\perp}\cdot v\,,
\end{array}
\]
where $\otimes$ denotes the vector tensor product.

The mathematical problem we are interested in reads: find the particle
temperature $u(t,x)$, solution of the evolution equation
\begin{gather} \label{P}
  (PH) \,\,\,
  \left\{
  \begin{array}{l}
    \partial_t u - {1 \over \eps} \nabla_{\parallel} \cdot (A_{\parallel}
    \nabla_{\parallel} u) -\nabla_{\perp} \cdot (A_{\perp} \nabla_{\perp}
    u)=0\,, \quad \textrm{in} \quad [0,T] \times \Omega\,,\\[3mm] 
    {1 \over \eps} n_{\parallel} \cdot (A_{\parallel}
    \nabla_{\parallel} u(t,\cdot)) + n_{\perp} \cdot (A_{\perp}
    \nabla_{\perp} u(t,\cdot))= g_N(t,\cdot)\,, \quad \textrm{on} \quad
          [0,T] \times \Gamma_{N}\,, \nonumber\\[3mm] 
          u (t,\cdot) = g_D(t,\cdot) \,, \quad \textrm{on} \quad [0,T] \times \Gamma_{D}\,, \nonumber\\[3mm]
          u(0,\cdot)=u^{0} (\cdot)\,, \quad \textrm{in} \quad \Omega\,.
  \end{array}
  \right.
  \nonumber
\end{gather}
%% The above problem describes the diffusion of an initial temperature
%% $u^{0}$ within the time interval $[0,T]$ and its outflow through the
%% boundaries. 
The diffusion coefficients $A_{\parallel}$ and $A_{\perp}$ are bounded
and of the same order of magnitude, satisfying
\begin{gather}
  0<A_{0}\leq A_{\parallel }(x)\leq A_{1}\,,\quad \text{for almost all}\,\,\,x\in
  \Omega ,
  % \label{eq:J48a1}
  \nonumber 
  \\[3mm]
  A_{0}||v||^{2}\leq v^{t}A_{\perp }(x)v\leq A_{1}||v||^{2}\,,\quad \forall
  v\in \mathbb{R}^{d}\,\,\,\text{and}\,\,\,\text{for almost all}\,\,\,x\in \Omega ,
  % \label{eq:J48a3}
  \nonumber 
\end{gather}%
with $0<A_{0}<A_{1}$ some positive constants. The anisotropy of the
problem is characterized by a parameter $\varepsilon $, which can be
very small and provoke substantial difficulties in the limit
$\varepsilon \to 0$.

Indeed, putting formally $\eps =0$ in (PH) leads to the following
reduced problem
\begin{gather}
  % \label{R}
  \nonumber 
  \left\{
  \begin{array}{l}
    - \nabla_{\parallel} \cdot (A_{\parallel} \nabla_{\parallel} u)
    =0\,, \quad \textrm{in} \quad [0,T] \times \Omega\,,
    \\[3mm]
    n_{\parallel} \cdot (A_{\parallel} \nabla_{\parallel} u(t,\cdot))
    = 0\,, \quad \textrm{on} \quad [0,T] \times \Gamma_{N}\,, \nonumber\\[3mm]
    u (t,\cdot) =g_D(t,\cdot) \,, \quad \textrm{on} \quad [0,T] \times \Gamma_{D}\,, \nonumber\\[3mm]
    u(0,\cdot)=u^{0} (\cdot)\,, \quad \textrm{in} \quad \Omega\,
  \end{array}
  \right. .
  \nonumber
\end{gather}
This problem may be ill-posed, depending on boundary conditions and
the anisotropy field $b$. For example, when some field lines of $b$
are closed in $\Omega $ the system would admit infinitely many
solutions as any function constant along the closed lines of $b$
(meaning $\nabla_{\parallel} u \equiv 0$) and satisfies the boundary
conditions solves the reduced problem. The same problem occurs when
the field lines are open but do not pass through a boundary supplied
with the Dirichlet conditions. This argument applies also to the case,
where periodic conditions are imposed on a part of the boundary. This
is the case in the numerical simulations related to the tokamak fusion
plasma, where computational domain is topologically equivalent to a
torus. Numerical discretization of the original (PH) problem in the
limit $\varepsilon \to 0$ can therefore lead to a very badly
conditioned linear systems. In fact, the condition number is
proportional to $1/\varepsilon $.

The goal of this work is to introduce a new numerical scheme which
permits to solve the original singular perturbation problem (PH)
independently of the value of $\varepsilon $ by the use of easy to
implement numerical tools without a need of adaptation of any kind to
the anisotropy strength/direction. The method presented herein
generalizes the numerical scheme introduced in \cite{LNN}. It removes
the limitations of the cited numerical scheme while keeping all the
advantages. It permits to solve the initial singular perturbation
problem (PH) accurately on a simple Cartesian grid, which does not
need to be specially tailored to suit the topology of the anisotropy
field lines $b$ and whose size does not depend on the anisotropy
strength.

The method is developed in the framework of Asymptotic-Preserving
schemes. That is to say, the method is consistent with a so-called
Limit problem as $\varepsilon \to 0$ and is stable independently of
the small parameter $\eps$. The initial singular perturbation problem
(PH) is reformulated in a suitable way. Introduction of an auxiliary
variable allows to remove the terms proportional to $1/\varepsilon $
from the equation. Resulting system is equivalent to the initial one
and is well-posed in the limit $\eps \rightarrow 0$. The modification
introduced here in allows to extend the applicability of the numerical
scheme proposed in \cite{LNN} to the settings, where the field $b$ may
contain closed lines without any loss of accuracy.

%%%%%%%%%%%%%%%%%%%%%%%%%%%%%%%%%%%%%%%%%%%%%%%
\section{Numerical method} \label{SEC3}
%%%%%%%%%%%%%%%%%%%%%%%%%%%%%%%%%%%%%%%%%%%%%%%

%%%%%%%%%%%%%%%%%%%%%%%%%%%%%%%%%%%%%%%%%%%%%%%
\subsection{Semi-discretization in space} \label{SEC32}
%%%%%%%%%%%%%%%%%%%%%%%%%%%%%%%%%%%%%%%%%%%%%%%

Let us write the variational formulation of the singular perturbation
problem (PH) : find $u(t,\cdot) \in \mathcal{V}:=H^1(\Omega)$ such
that
\begin{eqnarray}
  (P)\,\,\,
  &&\langle \partial _{t}u(t,\cdot ),\,v \rangle
  _{\mathcal{V}^{\ast },\mathcal{V}}+{1 \over \eps} \int_{\Omega
  }A_{\parallel}\nabla
  _{\parallel}u(t,\cdot )\cdot \nabla _{\parallel}v\, dx  \label{weak_bis} \\
  && \hspace{2cm} +\int_{\Omega }A_{\perp }\nabla _{\perp }u(t,\cdot )\cdot \nabla _{\perp
  }v\, dx
  - \int_{\Gamma_N } g_N(t,\cdot ) v
  =0,\quad \forall v \in \mathcal{V}  \notag\,
  \nonumber
\end{eqnarray}
for almost every $t \in (0,T)$. As already discussed in the previous
section, taking the formal limit of $\varepsilon \to 0$ leads to an
ill-posed problem:
\begin{equation*}
  \int_{\Omega }A_{\parallel}\nabla _{\parallel}u(t,\cdot )\cdot \nabla _{\parallel}v\, dx =0
\end{equation*}
with any function belonging to the vector space of functions constant
in the anisotropy direction:
\begin{gather}  
  {\cal G}:=\{p\in {\cal V}~/~\nabla _{\parallel}p=0 \text{ in }  \Omega \}
  \nonumber
\end{gather}
being a solution.

As proposed in \cite{LNN}, a correct Limit problem can be established
by seeking a solution in the subspace $\mathcal G$ instead of
$\mathcal V$. In this case the leading order term (containing the
parallel gradient) is eliminated from the equation and we are left
with the following Limit problem: find $u(t,\cdot) \in \mathcal{G}$
such that
\begin{gather}
  (L)\,\,\,
  \langle \partial _{t}u(t,\cdot ),\,v \rangle
  _{\mathcal{V}^{\ast },\mathcal{V}}
  +\int_{\Omega }A_{\perp }\nabla _{\perp }u(t,\cdot )\cdot \nabla _{\perp
  }v\, dx
  - \int_{\Gamma_N } g_N(t,\cdot ) v
  =0,\quad \forall v \in \mathcal{G}  
  \nonumber
\end{gather}
for almost every $t \in (0,T)$.

An Asymptotic Preserving scheme designed for the problem (P) should
give accurate results and be stable independently of $\varepsilon $,
even in the limit $\varepsilon =0$. In particular, the condition
number associated with corresponding linear system should not depend
on $\varepsilon $. That is to say, the correct Limit problem (L) has
to be ``hardcoded'' in certain sense into the equations. Let us
briefly recall the Asymptotic-Preserving scheme introduced in
\cite{LNN}. The method relies on the auxiliary variable $q$ introduced
by the relation $\eps\nabla_{\parallel}q=\nabla_{\parallel}u$ in
$\Omega$. This trick permits to get rid of the terms of order
$O(1/\eps)$ from the variational formulation. The uniqueness of $q$ is
provided by setting $q=0$ on a part of the boundary $\Gamma_{in}$
defined by
\begin{equation*}
  \Gamma_{in}:= \{ x \in \Gamma \,\, / \,\, b(x) \cdot n(x) <0 \}\,,
\end{equation*}
{\it i.e.} the part of the boundary, where the field lines enter the
domain. The following reformulated problem, called in the sequel the
Asymptotic-Preserving reformulation (AP-problem) is proposed: find
$(u(t,\cdot),q(t,\cdot))\in {\cal V} \times {\cal L_{in}}$, solution
of
\begin{equation}
  (AP)\,\,\, 
  \left\{ 
  \begin{array}{l}
    \ds \langle\frac{\partial u}{\partial t},\, v \rangle_{\mathcal{V}^{\ast },\mathcal{V}} + \int_{\Omega }(
    A_{\perp }\nabla _{\perp }u)\cdot \nabla _{\perp
    }v\,dx+\int_{\Omega } A_{\parallel} \nabla _{\parallel}q\cdot
    \nabla_{\parallel}v\,dx
    - \int_{\Gamma_{N}} g_Nv\,ds
    % =\int_{\Omega }fv\,dx, & \quad \forall v\in {\cal V}
    =0, \\[1mm]
    \hspace{12cm}\forall v\in {\cal V}
    \\[1mm] \ds \int_{\Omega }A_{\parallel} \nabla_{\parallel}u\cdot
    \nabla_{\parallel}w\,dx-\varepsilon \int_{\Omega } A_{\parallel} \nabla
    _{\parallel}q\cdot \nabla _{\parallel}w\,dx=0,\quad \forall w\in {\cal
      L_{in}}\,,
  \end{array}
  \right.  \label{Pa}
\end{equation}
where 
\begin{gather}
  {\cal L_{in}}:=\{q\in L^{2}(\Omega )~/~\nabla _{\parallel}q\in
  L^{2}(\Omega )\text{ and }  q|_{\Gamma _{in} }=0\}. 
  % \label{eq:Jobb}
\nonumber 
\end{gather}
The AP-problem is equivalent for fixed $\eps>0$ to the original
P-problem (\ref{weak_bis}). Moreover, putting formally $\eps =0$ in
(AP) leads to a well-posed problem
\begin{equation*}
  (L')\,\,\, 
  \left\{ 
  \begin{array}{l}
    \ds \langle\frac{\partial u}{\partial t},\, v \rangle_{\mathcal{V}^{\ast },\mathcal{V}} + \int_{\Omega }(
    A_{\perp }\nabla _{\perp }u)\cdot \nabla _{\perp
    }v\,dx+\int_{\Omega } A_{\parallel} \nabla _{\parallel}q\cdot
    \nabla_{\parallel}v\,dx
    - \int_{\Gamma_{N}} g_Nv\,ds
    % =\int_{\Omega }fv\,dx, & \quad \forall v\in {\cal V}
    =0,\\[1mm]
    \hspace{12cm} \forall v\in {\cal V}
    \\[1mm] \ds \int_{\Omega }A_{\parallel} \nabla_{\parallel}u\cdot
    \nabla_{\parallel}w\,dx=0,\quad \forall w\in {\cal
      L_{in}}\,,
  \end{array}
  \right.
  % \label{LP}
\end{equation*}
which is equivalent to the correct Limit problem (L). In this case,
the auxiliary variable $q$ acts as a Lagrange multiplier forcing $u$
to be constant along $b$.

The drawback of this method is the choice of the space for the
auxiliary variable. Imposing $q|_{\Gamma_{in}} = 0$ provides
uniqueness of a solution but limits the application of the scheme to
the case where all field lines are open. Indeed, fixing a value of $q$
on the inflow boundary does not provide uniqueness of $q$ on field
lines which does not intersect with the inflow boundary ({\it i.e.} on
closed field lines). In order to overcome this restriction we propose
a new approach based on penalty stabilization rather than on fixing
the value of $q$ on one of the boundaries. The modification of the
second equation of the AP scheme (\ref{Pa}) consists of an
introduction of a penalty term --- a mass matrix $\int_{\Omega} qw$
multiplied by a stabilization constant. A suitable choice of this
constant permits to conserve an accuracy of the scheme.

The new method (APS-scheme) reads: find $(u(t,\cdot),q(t,\cdot))\in
{\cal V} \times {\cal L}$, solution of
\begin{equation}
  (APS)\,\,\, 
  \left\{ 
  \begin{array}{l}
    \ds \langle\frac{\partial u}{\partial t},\, v \rangle_{\mathcal{V}^{\ast },\mathcal{V}} + \int_{\Omega }(
    A_{\perp }\nabla _{\perp }u)\cdot \nabla _{\perp
    }v\,dx+\int_{\Omega } A_{\parallel} \nabla _{\parallel}q\cdot
    \nabla_{\parallel}v\,dx
    - \int_{\Gamma_{N}} g_Nv\,ds
    % =\int_{\Omega }fv\,dx, & \quad \forall v\in {\cal V}
    =0, \\[1mm]
    \hspace{12cm}\forall v\in {\cal V}
    \\[1mm] \ds \int_{\Omega }A_{\parallel} \nabla_{\parallel}u\cdot
    \nabla_{\parallel}w\,dx-\varepsilon \int_{\Omega } A_{\parallel} \nabla
    _{\parallel}q\cdot \nabla _{\parallel}w\,dx - \alpha 
    \int_{\Omega} qw=0,\quad \forall w\in {\cal
      L\,},
  \end{array}
  \right.  \label{APS}
\end{equation}
where 
\begin{gather}
  {\cal L}:=\{q\in L^{2}(\Omega )~/~\nabla _{\parallel}q\in
  L^{2}(\Omega )\} 
  % \label{eq:Jobb}
\nonumber 
\end{gather}
and $\alpha$ is a positive stabilization constant. We postpone a
theoretical justification of the proposed method to the forthcoming
paper. A similar method with stabilization by a diffusion matrix
instead of a mass matrix was recently proposed in the context of a
{\it a posteriori} error indicator and mesh adaptation for strongly
anisotropic elliptic equations in \cite{N}.

Let us now choose a polygonalization of the domain $\Omega $ with
polygons of the diameter approximately equal to $h$ and introduce the
finite element spaces ${\cal V}_{h} \subset {\mathcal V}$ and ${\cal
  L}_{h}\subset {\mathcal L}$. The finite element discretization of
(\ref{APS}) writes then: find $(u_{h},q_{h}) \in {\cal V}_{h} \times
{\cal L}_h$ such that

\begin{equation}
  (APS)_h\,\,\, 
  \left\{ 
  \begin{array}{l}
    \ds \int_{\Omega } \frac{\partial u_h}{\partial t}v_h \, dx + \int_{\Omega }(
    A_{\perp }\nabla _{\perp }u_h)\cdot \nabla _{\perp
    }v_h\,dx+\int_{\Omega } A_{\parallel} \nabla _{\parallel}q_h\cdot
    \nabla_{\parallel}v_h\,dx
    - \int_{\Gamma_{N}} g_Nv_h\,ds
    %% +\gamma \int_{\Gamma_{\perp}} u_hv_h\,ds
    %      =\int_{\Omega }fv\,dx, & \quad \forall v\in {\cal V}
    =0, \\[1mm]
    \hspace{12cm} \forall v_h\in {\cal V}_h
    \\[1mm] \ds \int_{\Omega }A_{\parallel} \nabla_{\parallel}u_h\cdot
    \nabla_{\parallel}w_h\,dx-\varepsilon \int_{\Omega } A_{\parallel} \nabla
    _{\parallel}q_h\cdot \nabla _{\parallel}w_h\,dx - h^{k+1}  \int_{\Omega} qw =0, \quad \forall w\in {\cal
      L}_h\,.
  \end{array}
  \right.
  \label{APSh}
\end{equation}
Remark that in order to ensure convergence rate in $L^2$-norm we have
put $\alpha = h^{k+1}$, where $k$ is the order of finite element
method.

%%%%%%%%%%%%%%%%%%%%%%%%%%%%%%%%%%%%%%%%%%%%%%%
\subsection{Semi-discretization in time} \label{SEC31}
%%%%%%%%%%%%%%%%%%%%%%%%%%%%%%%%%%%%%%%%%%%%%%

One should be extremely careful when discretizing in time the
(\ref{APSh}) scheme as not all numerical schemes conserve the AP
property. In fact a chosen method should be L-stable. This is the case
for a standard first order implicit Euler scheme. If however a higher
order in time method is desired then certain simple schemes are
excluded. For example a Crank-Nicolson discretization is only A-stable
and not L-stable. The obtained system would not be asymptotic
preserving giving reliable results only under certain
assumptions. Namely, $\eps$ should be close to one, time step should
be of the order of $\eps$ or the initial value $u^0$ should already be
a solution to the stationary equation, {\it i.e.} the parallel
gradient of $u^0$ should be proportional to $\eps$. This is why we
choose in this work to present both a first order implicit Euler
scheme and a second order, L-stable Runge-Kutta method.

\subsubsection{Implicit Euler scheme} 
\quad\\

Let us introduce the standard bilinear forms
\begin{gather}
  (\Theta,\chi):= \int_{\Omega} \Theta \chi \, dx 
  \,,
  % \label{eq:Jmbb0}
  \nonumber 
  \\
  a_\parallel(\Theta,\chi):= \int_{\Omega} A_\parallel \nabla_\parallel \Theta \cdot
  \nabla_\parallel \chi \, dx 
  \,, \quad  \quad a_{\perp}(\Theta,\chi):=
  \int_{\Omega} A_{\perp} \nabla_{\perp} \Theta \cdot \nabla_{\perp}
  \chi \, dx \,,
  % \label{eq:Jmbb3}
  \nonumber 
\end{gather}
and write the first order, implicit Euler method in more compact
notation: Find $(u_h^{n+1},q_h^{n+1}) \in {\mathcal V}_h \times
{\mathcal L}_h$, solution of
\begin{gather}
  (E_{APS}) \,\,\, 
  \left\{
  \begin{array}{l}
    (u_h^{n+1},v_h) + 
    \tau \left(a_\perp (u_h^{n+1},v_h) + a_\parallel (q_h^{n+1},v_h) 
    - \int_{\Gamma_{N}} g_N(t^{N+1})v_h\,ds
    %% + \gamma \int_{\Gamma_{\perp}} u_h^{n+1}v_h \, ds
    \right)
    = (u_h^{n},v_h)
    \nonumber\\[4mm]
    a_{\parallel} (u_h^{n+1},w_h) - \eps a_\parallel
    (q_h^{n+1},w_h)
    -h^{k+1} (q_h^{n+1},w_h)
    = 0\,,
  \end{array}
  \right.
  %  \label{eq:Jpcb}
  \nonumber 
  .
\end{gather}

\subsubsection{L-stable Runge-Kutta method}

Any $s$-stage Runge-Kutta method applied to following problem
\begin{gather}
  \frac{\partial u}{\partial t} = Lu + f(t)\,,
  % \label{eq:J6cb}
  \nonumber 
\end{gather}
is conveniently defined using a Butcher diagram:
\begin{gather}
  \begin{array}{c|ccc}
    c_1 & a_{11} & \cdots & a_{1s} \nonumber\\
    \vdots & \vdots & & \vdots \nonumber\\
    c_s & a_{s1} & \cdots & a_{ss}\nonumber\\
    \hline
    & b_1 & \cdots & b_s
  \end{array}.
  \nonumber
\end{gather}
The method reads: for given $u^n$, an approximation of $u(t_n)$,
the $u^{n+1}$ is a linear combination
the method
\begin{gather}
  u^{n+1} = u^n + \tau \sum_{j=1}^{s} b_j u_j
  \nonumber,
\end{gather}
where $u_i$ are solutions to the following problems:
\begin{gather}
  u_i = u^{n} + \tau \sum_{j=1}^{s} a_{ij} (Lu_j + f(t+ c_j\tau)) ,
  % \label{eq:J8cb}
  \nonumber 
  .
\end{gather}
Moreover, if $b_j = a_{sj}$ for $j=1,\ldots, s$ than $u^{n+1}$ is
equal to the last stage of the method: $u^{n+1} = u_s$.

In order to obtain a second order accurate in time scheme, we choose
to implement a two stage Diagonally Implicit Runge-Kutta (DIRK) second
order scheme. The scheme is developed according to the following
Butcher's diagram:
\begin{gather}
  \begin{array}{c|cc}
    \lambda & \lambda & 0 \nonumber\\
    1 & 1-\lambda & \lambda \nonumber\\
    \hline
    & 1-\lambda & \lambda 
  \end{array}
  \label{eq:Jrcb}.
\end{gather}
The method is known to be L-stable for $\lambda = 1-\frac{1}{\sqrt{2}}$. 

The second order AP-scheme writes: find $(u_h^{n+1},q_h^{n+1}) \in
{\mathcal V}_h \times {\mathcal L}_h$, solution of
\begin{gather}
  (RK_{APS}) \,\,\, 
  \left.
  \begin{array}{l}
    \left\{
    \begin{array}{l}
      (u_{1,h}^{n+1},v_h)+ 
      \tau \lambda \left(a_\perp (u_{1,h}^{n+1},v_h) 
      + a_\parallel (q_{1,h}^{n+1},v_h) \right)    
      - \int_{\Gamma_{N}} g_N(t^N + \lambda \tau )v_h\,ds
      \nonumber\\[3mm]
      \qquad \qquad \qquad = (u_{h}^{n},v_h)    
      \nonumber\\[4mm]
      a_{\parallel}
      \left(u_{1,h}^{n+1},w_h\right)
      - \eps a_\parallel (q_{1,h}^{n+1},w_h) - h^{k+1} (q_{1,h}^{n+1},w_h)= 0
    \end{array}
    \right.
    \nonumber\\
    \nonumber\\
    \left\{
    \begin{array}{l}
      (u_{2,h}^{n+1},v_h)+ 
      \tau \lambda \left(a_\perp (u_{2,h}^{n+1},v_h) 
      + a_\parallel (q_{2,h}^{n+1},v_h) \right)      
      - \int_{\Gamma_{N}} g_N(t^{N}+\tau )v_h\,ds 
      \nonumber\\[3mm]
      \qquad \qquad \qquad 
      = (u_{h}^{n},v_h) + {1-\lambda \over \lambda } \left( u_{1,h}^{n+1}-u_{h}^{n},v_h\right)
      \nonumber\\[4mm]
      a_{\parallel}
      \left(u_{2,h}^{n+1},w_h\right)
      - \eps a_\parallel (q_{2,h}^{n+1},w_h) - h^{k+1} (q_{2,h}^{n+1},w_h)= 0
    \end{array}
    \right.
    \nonumber\\
    \nonumber\\
    u_h^{n+1} = u_{2,h}^{n+1}\,, \qquad q_h^{n+1} = q_{2,h}^{n+1}\,,
    \end{array}
  \right.
  % \label{eq:Jscb}
  \nonumber 
\end{gather}
with $u_{1,h}^{n+1}$ and $u_{2,h}^{n+1}$ denote the solutions of the
first and the second stage of the Runge-Kutta method.

%%%%%%%%%%%%%%%%%%%%%%%%%%%%
\section{Numerical results}\label{sec:numres}
%%%%%%%%%%%%%%%%%%%%%%%%%%%%

We are now ready to test the proposed schemes. Let us perform in the
beginning numerical experiments for $b$ having all field lines
open. This setting allows us to compare implicit Euler-APS and
DIRK-APS methods with their non penalty stabilized versions (Euler-AP
and DIRK-AP) proposed in \cite{LNN}. The Euler-AP and DIRK-AP schemes
are obtained by discretization of the AP-problem and differ from
Euler-APS and DIRK-APS in two details. Stabilization terms are absent
and the vector space $\cal L _h$ is replaced by $\cal L_{in,h}$.  The
methods are also compared to a standard implicit Euler discretization
of the initial singular perturbation problem (PH), given by
\begin{gather}
  (P)_{h\tau } \quad
  (u_h^{n+1},v_h) + 
  \tau \left(a_\perp (u_h^{n+1},v_h) + \frac{1}{\eps}a_{\parallel nl} (u_h^{n},u_h^{n+1},v_h) 
    - \int_{\Gamma_{N}} g_N(t^{N+1})v_h \, ds \right)
  = 
  (u_h^{n},v_h)  
  % \label{eq:Jtcb}
  \nonumber 
  .
\end{gather}

Finally, the DIRK-APS scheme is applied to the case of magnetic
islands, where some field lines of $b$ are closed. But let us first
introduce a finite element space that we intend to use in all
numerical experiments.

%%%%%%%%%%%%%%%%%%%%%%%
\subsection{Discretization} \label{Discr}
%%%%%%%%%%%%%%%%%%%%%%%

Let us consider a 2D square computational domain $\Omega = [0,1]\times
[0,1]$. We choose to perform all simulations on structured meshes. Let
us introduce the Cartesian, homogeneous grid
\begin{gather}
  x_i = i / N_x \;\; , \;\; 0 \leq i \leq N_x \,, \quad
  y_j = j / N_y \;\; , \;\; 0 \leq j \leq N_y
%  \label{eq:Jp8a}
\nonumber 
,
\end{gather}
with $N_x$ and $N_y$ being positive even constants, corresponding to
the number of discretization intervals in the $x$-
resp. $y$-direction. The corresponding mesh-sizes are denoted by $h_x
>0$ resp. $h_y >0$. We choose a standard $\mathbb Q_2$ finite element
method ($\mathbb Q_2$-FEM), based on the following quadratic base
functions

\begin{gather}
  \theta _{x_i}=
  \left\{
    \begin{array}{ll}
      \frac{(x-x_{i-2})(x-x_{i-1})}{2h_x^{2}} & x\in [x_{i-2},x_{i}],\nonumber\\
      \frac{(x_{i+2}-x)(x_{i+1}-x)}{2h_x^{2}} & x\in [x_{i},x_{i+2}],\nonumber\\
      0 & \text{else}
    \end{array}
  \right.\,, \quad 
  \theta _{y_j} =
  \left\{
    \begin{array}{ll}
      \frac{(y-y_{j-2})(y-y_{j-1})}{2h_y^{2}} & y\in [y_{j-2},y_{j}],\nonumber\\
      \frac{(y_{j+2}-y)(y_{j+1}-y)}{2h_y^{2}} & y\in [y_{j},y_{j+2}],\nonumber\\
      0 & \text{else}
    \end{array}
  \right.
%  \label{eq:Js8a1}
\nonumber 
\end{gather}
for even $i,j$  and
\begin{gather}
  \theta _{x_i}=
  \left\{
    \begin{array}{ll}
      \frac{(x_{i+1}-x)(x-x_{i-1})}{h_x^{2}} & x\in [x_{i-1},x_{i+1}],\nonumber\\
      0 & \text{else}
    \end{array}
  \right.\,, \quad 
  \theta _{y_j} =
  \left\{
    \begin{array}{ll}
      \frac{(y_{j+1}-y)(y-y_{j-1})}{h_y^{2}} & y\in [y_{j-1},y_{j+1}],\nonumber\\
      0 & \text{else}
    \end{array}
  \right.
%  \label{eq:Js8a2}
\nonumber 
\end{gather}
for odd $i,j$, we define the space
$$
\cal W_h := \{ v_h = \sum_{i,j} v_{ij}\,  \theta_{x_i} (x)\,  \theta_{y_j}(y)\}\,.
$$
The spaces ${\cal V}_h$ and ${\cal L}_h$ are then defined by
\begin{equation*}
  {\cal V}_{h}={\cal L}_{h}={\cal W}_{h}\,, \quad {\cal L}_{in,h}=\{q_{h}\in {\cal
    V}_{h}\,, \text{ such that\,\,\,} q_{h}|_{\Gamma _{in}}=0\}. 
\end{equation*}
The matrix elements are computed using the 2D Gauss quadrature
formula, with 3 points in the $x$ and $y$ direction:
\begin{gather}
  \int_{-1}^{1}\int_{-1}^{1}f (\xi,\eta)\, d\xi\, d\eta =
  \sum_{i,j=-1}^{1} \omega _{i}\omega _{j} f (\xi_i,\eta_j)\,,
%  \label{eq:Jk9a}
\nonumber 
\end{gather}
where $\xi_0=\eta_0=0$, $\xi_{\pm 1}=\eta_{\pm 1}=\pm\sqrt
{\frac{3}{5}}$, $\omega _0 = 8/9$ and $\omega _{\pm 1} = 5/9$, which
is exact for polynomials of degree 5. Linear systems obtained for all
methods in these numerical experiments are solved using a LU
decomposition, implemented by the MUMPS library.

%%%%%%%%%%%%%%%%%%%%%%%%%%%%
\subsection{Numerical tests}
%%%%%%%%%%%%%%%%%%%%%%%%%%%%

\subsubsection{Known analytical solution}

Let the computational domain $\Omega $ be supplied with the boundaries
$\Gamma _N = \{ (0,x) \cup (1,x) | x \in [0,1] \}$ and $\Gamma _D = \{
(x,0) \cup (x,1) | x \in [0,1] \}$ with the boundary conditions $g_N =
g_D =0$.

Let us now construct a numerical test case with a known
solution. Finding an analytical solution for an arbitrary $b$-field is
extremely difficulty. We prefer rather to act as in the previous
papers \cite{DDLNN,LNN}, where we started from selecting a suitable
limit solution $u_0$ and than defined anisotropy direction as its
isolines. 
Let
\begin{gather}
  u_{0} = \sin \left(\pi y +\alpha (y^2-y)\cos (\pi x) \right) e^{-t}
  % \label{eq:J79a}
  \nonumber 
  ,
\end{gather}
where $\alpha $ is an amplitude of variations of $b$. For $\alpha =0$,
the limit solution is constant in the $X$-direction and represents a
solution for $b=(1,0)^T$. For $\alpha \neq 0$ the anisotropy
direction is determined by the following implication
\begin{gather}
  \nabla_{\parallel} u_{0} = 0 \quad \Rightarrow \quad
  b_x \frac{\partial u_{0}}{\partial x} +
  b_y \frac{\partial u_{0}}{\partial y} = 0\,,
%  \label{eq:J89a}
\nonumber 
\end{gather}
and therefore $b$ can be defined for example as 
\begin{gather}
  b = \frac{B}{|B|}\, , \quad
  B =
  \left(
    \begin{array}{c}
      \alpha  (2y-1) \cos (\pi x) + \pi \nonumber\\
      \pi \alpha  (y^2-y) \sin (\pi x)
    \end{array}
  \right)
  % \label{eq:J99a}
  \nonumber 
  \,\quad. 
\end{gather}
We set $\alpha =1$ in all our simulations so that the direction of the
anisotropy is variable in the computational domain.  Note that we have
$B \neq 0$ in the computational domain.

Finally, a perturbation proportional to $\varepsilon $ is added to
$u^0$ and an a function $u $ is obtained. In this way we ensure that
$u$ converges, as $\varepsilon \rightarrow 0$, to the limit solution
$u_{0}$. For example
\begin{gather}
  u = \sin \left(\pi y +\alpha (y^2-y)\cos (\pi
  x) \right) e^{-t} + \varepsilon \cos \left( 2\pi x\right) \sin \left(\pi
  y \right) e^{-t}
  \label{eq:Jc0a}.
\end{gather}

We set the initial condition $u^0$ to be equal to $u(t=0)$, with $u$
defined by (\ref{eq:Jc0a}) and add a suitable force term to the right
hand side of the numerical schemes so that $u$ is an exact solution to
the problem. In this setting we expect all Asymptotic-Preserving
methods to converge in the optimal rate, independently on $\eps$.

\begin{table}
  \centering
  \begin{tabular}{|c||c|c|c|c|c|}
    \hline%\rule{0pt}{2.5ex}
    \multirow{2}{*}{$h$}  
    & \multicolumn{5}{|c|}{\rule{0pt}{2.5ex}$L^2$-error\qquad $\eps=1$} \\
    \cline{2-6} 
    &\rule{0pt}{2.5ex}
    P& $E_{AP}$& $E_{APS}$& $RK_{AP}$& $RK_{APS}$\\
    \hline
    \hline\rule{0pt}{2.5ex}
    0.1 &
    $5.6\times 10^{-3}$ &
    $5.6\times 10^{-3}$ &
    $5.6\times 10^{-3}$ &
    $5.6\times 10^{-3}$ &
    $5.6\times 10^{-3}$ 
    \\
    \hline\rule{0pt}{2.5ex}
    0.05 &
    $7.1\times 10^{-4}$ &
    $7.1\times 10^{-4}$ &
    $7.1\times 10^{-4}$ &
    $7.1\times 10^{-4}$ &
    $7.1\times 10^{-4}$ 
    \\
    \hline\rule{0pt}{2.5ex}
    0.025 &
    $8.9\times 10^{-5}$ &
    $8.9\times 10^{-5}$ &
    $8.9\times 10^{-5}$ &
    $8.9\times 10^{-5}$ &
    $8.9\times 10^{-5}$ 
    \\
    \hline\rule{0pt}{2.5ex}
    0.0125 &
    $1.11\times 10^{-5}$ &
    $1.11\times 10^{-5}$ &
    $1.11\times 10^{-5}$ &
    $1.11\times 10^{-5}$ &
    $1.11\times 10^{-5}$ 
    \\
    \hline\rule{0pt}{2.5ex}
    0.00625 &
    $1.39\times 10^{-6}$ &
    $1.39\times 10^{-6}$ &
    $1.39\times 10^{-6}$ &
    $1.39\times 10^{-6}$ &
    $1.39\times 10^{-6}$ 
    \\
    \hline\rule{0pt}{2.5ex}
    0.003125 &
    $1.74\times 10^{-7}$ &
    $1.74\times 10^{-7}$ &
    $1.74\times 10^{-7}$ &
    $1.74\times 10^{-7}$ &
    $1.74\times 10^{-7}$ 
    \\
    \hline
  \end{tabular}
  \begin{tabular}{|c||c|c|c|c|c|}
    \hline%\rule{0pt}{2.5ex}
    \multirow{2}{*}{$h$}  
    & \multicolumn{5}{|c|}{\rule{0pt}{2.5ex}$L^2$-error\qquad $\eps=10^{-20}$} \\
    \cline{2-6} 
    &\rule{0pt}{2.5ex}
    P& $E_{AP}$& $E_{APS}$& $RK_{AP}$& $RK_{APS}$\\
    \hline
    \hline\rule{0pt}{2.5ex}
    0.1 &
    $6.9\times 10^{-1}$ &
    $1.62\times 10^{-3}$ &
    $1.66\times 10^{-3}$ &
    $1.62\times 10^{-3}$ &
    $1.66\times 10^{-3}$ 
    \\
    \hline\rule{0pt}{2.5ex}
    0.05 &
    $6.9\times 10^{-1}$ &
    $2.20\times 10^{-4}$ &
    $2.37\times 10^{-4}$ &
    $2.20\times 10^{-4}$ &
    $2.37\times 10^{-4}$
    \\
    \hline\rule{0pt}{2.5ex}
    0.025 &
    $6.9\times 10^{-1}$ &
    $2.77\times 10^{-5}$ &
    $2.93\times 10^{-5}$ &
    $2.77\times 10^{-5}$ &
    $2.93\times 10^{-5}$
    \\
    \hline\rule{0pt}{2.5ex}
    0.0125 &
    $6.9\times 10^{-1}$ &
    $3.43\times 10^{-6}$ &
    $3.58\times 10^{-6}$ &
    $3.43\times 10^{-6}$ &
    $3.58\times 10^{-6}$
    \\
    \hline\rule{0pt}{2.5ex}
    0.00625 &
    $6.9\times 10^{-1}$ &
    $4.2\times 10^{-7}$ &
    $4.4\times 10^{-7}$ &
    $4.2\times 10^{-7}$ &
    $4.4\times 10^{-7}$
    \\
    \hline\rule{0pt}{2.5ex}
    0.003125 &
    $6.9\times 10^{-1}$ &
    $5.3\times 10^{-8}$ &
    $5.5\times 10^{-8}$ &
    $5.3\times 10^{-8}$ &
    $5.5\times 10^{-8}$
    \\
    \hline
  \end{tabular}
 
  \caption{The absolute error of $u$ in the $L^{2}$-norm for different
    mesh sizes and $\eps =1$ or $\eps =10^{-20}$, using the singular
    perturbation scheme (P) and the two proposed AP-schemes for a time
    step of $\tau = 10^{-6}$ and at instant $t=10^{-4}$.  }
  \label{tab:conv_h1}
\end{table}

In order to validate our method and confirm our expectations we verify
first the convergence in regard of the space discretization. We choose
a time step small enough so that the space discretization error is
much bigger than the time discretization error. Then, we perform
numerical simulations for 100 time steps for different mesh sizes for
all Asymptotic-Preserving schemes and the implicit Euler
discretization of the (P) problem. Numerical errors are given in Table
\ref{tab:conv_h1} and Figures \ref{fig:conv_h1} and
\ref{fig:conv_h2}. The third order space convergence in the $L_2$-norm
for large values of $\eps$ is achieved (as expected) by all
methods. Stabilization procedure does not alter the accuracy for weak
anisotropy. For small values of $\eps$ only the Asymptotic Preserving
schemes give good numerical solutions. In this case the stabilization
term decreases slightly the precision of the results (by $2.9 - 4.6
\%$ compared to non stabilized schemes), keeping however the order of
convergence.

\def\xxxa{0.45\textwidth}
\begin{figure}[!ht]
  \centering
  \subfigure[$\eps = 1$]
  {\includegraphics[angle=0,width=\xxxa]{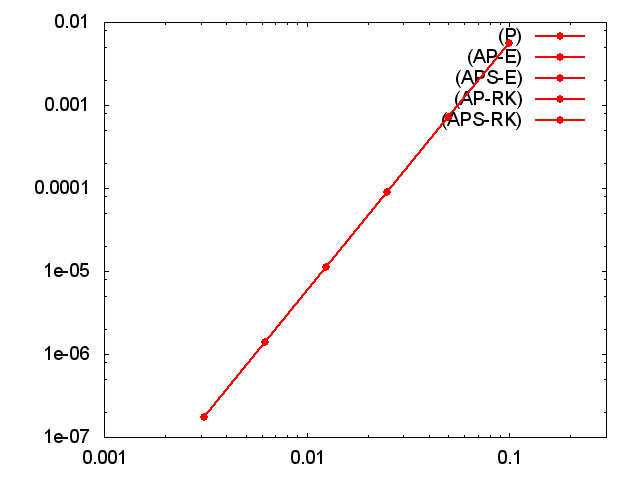}}
  \subfigure[$\eps = 10^{-20}$]
  {\includegraphics[angle=0,width=\xxxa]{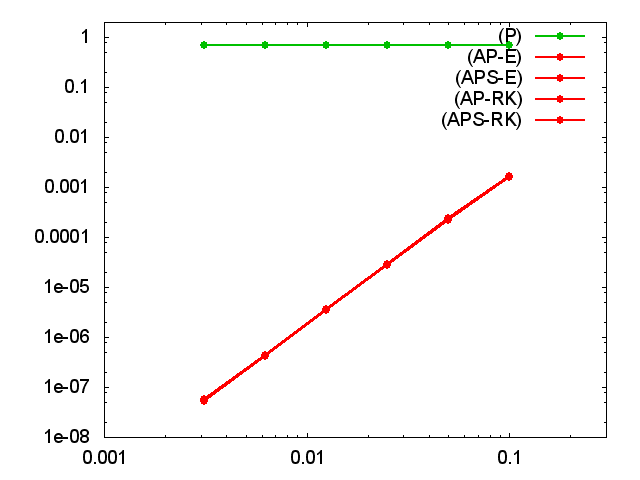}}
  \caption{Relative $L^{2}$-errors between the exact solution
    $u^{\varepsilon }$ and the computed solution for the standard
    scheme (P), the Euler-AP method ($E_{AP}$) the stabilized
    Euler-APS method ($E_{APS}$), the DIRK-AP scheme ($RK_{AP}$) and
    the stabilized DIRK-APS scheme ($RK_{APS}$) as a function of $h$,
    for $\eps=1$ resp. $\eps=10^{-20}$ and the time step
    $\tau=10^{-6}$. Observe that for $\eps=1$ all schemes give the
    same precision, for $\eps = 10^{-20}$ the standard scheme does not
    work while all AP schemes give comparable accuracy.}
  \label{fig:conv_h1}
\end{figure}

\def\xxxa{0.45\textwidth}
\begin{figure}[!ht]
  \centering
  \subfigure[$h = 0.1$]
  {\includegraphics[angle=0,width=\xxxa]{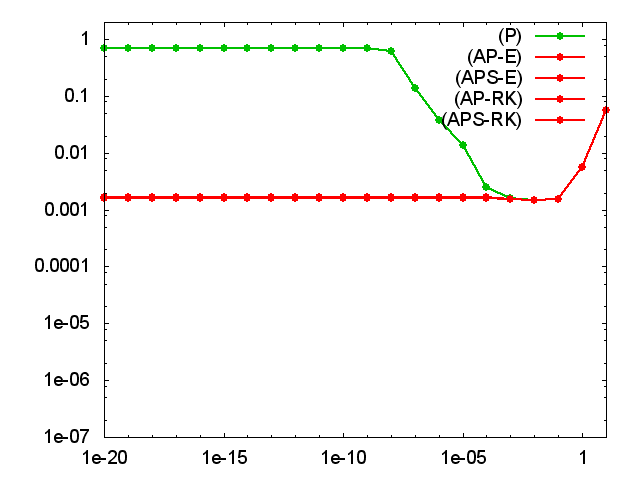}}
  \subfigure[$h = 0.00625$]
  {\includegraphics[angle=0,width=\xxxa]{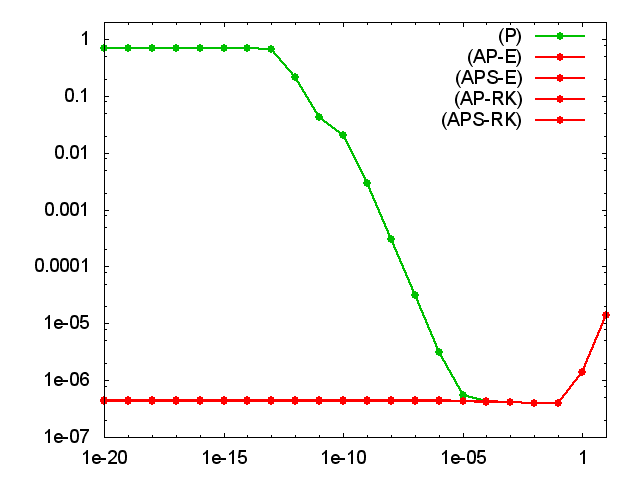}}
  \caption{Relative $L^{2}$-errors between the exact
    solution $u^{\varepsilon }$ and the computed solution for the standard
    scheme (P), the Euler-AP method ($E_{AP}$) the stabilized
    Euler-APS method ($E_{APS}$), the DIRK-AP scheme ($RK_{AP}$) and
    the stabilized DIRK-APS scheme ($RK_{APS}$)
    as a function of
    $\eps$, for $h=0.1$ resp. $h=0.00625$ and the time step
    $\tau=10^{-6}$. Observe that the standard scheme is accurate only
    for small values of anisotropy (large $\eps$)
    while all AP schemes give comparable accuracy in the whole range
    of $\eps$.}
  \label{fig:conv_h2}
\end{figure}

Next, the time convergence of the methods is tested. In this case the
mesh size is chosen in such a way that the space discretization does
not influence the numerical precision (at least for large time steps).
We perform simulations for different time steps and for a fixed final
time ($t=0.1$). The results are presented in Table \ref{tab:conv_t}
and Figure \ref{fig:conv_t2}. The time discretization convergence
order is confirmed. The $(RK_{AP})$ and $(RK_{APS})$ schemes are, as
expected, of second order in time as long as the error due to the
space discretization is smaller than the error induced by the time
discretization. The $(E_{AP})$ and $(E_{APS})$ schemes are of first
order for all values of the anisotropy parameter and the standard
(P)-scheme works well and is of first order where the anistropy is
weak ($\eps$ close to one). Note that the stabilization procedure
does not influence the accuracy of the solution if the time
discretization error is bigger than the space discretization
error. One can also observe that reasonably accurate results are
obtained even with relatively large time steps, especially for the
Runge-Kutta schemes. This property would allow orders of magnitude
gains in integration time for a given physics application with respect
to non-AP schemes.

\begin{table}
  \centering
  \begin{tabular}{|c||c|c|c|c|c|}
    \hline%\rule{0pt}{2.5ex}
    \multirow{2}{*}{$\tau$}  
    & \multicolumn{5}{|c|}{\rule{0pt}{2.5ex}$L^2$-error\qquad $\eps=1$} \\
    \cline{2-6} 
    &\rule{0pt}{2.5ex}
    P& $E_{AP}$& $E_{APS}$& $RK_{AP}$& $RK_{APS}$\\
    \hline
    \hline\rule{0pt}{2.5ex}
    0.1 &
    $1.32\times 10^{-3} $&
    $ 1.32\times 10^{-3} $&
    $ 1.32\times 10^{-3} $&
    $ 8.4\times 10^{-5} $&
    $ 8.4\times 10^{-5} $
    \\
    \hline\rule{0pt}{2.5ex}
    0.05 &
    $7.2\times 10^{-4} $&
    $7.2\times 10^{-4} $&
    $7.2\times 10^{-4} $&
    $2.43\times 10^{-5} $&
    $2.43\times 10^{-5} $
     \\
    \hline\rule{0pt}{2.5ex}
    0.025 &
    $3.55\times 10^{-4} $&
    $3.55\times 10^{-4} $&
    $3.55\times 10^{-4} $&
    $6.2\times 10^{-6} $&
    $6.2\times 10^{-6} $
    \\
    \hline\rule{0pt}{2.5ex}
    0.0125 &
    $1.75\times 10^{-4} $&
    $1.75\times 10^{-4} $&
    $1.75\times 10^{-4} $&
    $1.59\times 10^{-6} $&
    $1.59\times 10^{-6} $
    \\
    \hline\rule{0pt}{2.5ex}
    0.00625 &
    $8.8\times 10^{-5} $&
    $8.8\times 10^{-5} $&
    $8.8\times 10^{-5} $&
    $4.8\times 10^{-7} $&
    $4.8\times 10^{-7} $
    \\
    \hline\rule{0pt}{2.5ex}
    0.003125 &
    $4.4\times 10^{-3} $&
    $4.4\times 10^{-3} $&
    $4.4\times 10^{-3} $&
    $2.82\times 10^{-5} $&
    $2.82\times 10^{-5} $
    \\
    \hline\rule{0pt}{2.5ex}
    0.0015625 &
    $2.20\times 10^{-5} $&
    $2.20\times 10^{-5} $&
    $2.20\times 10^{-5} $&
    $2.64\times 10^{-7} $&
    $2.64\times 10^{-7} $
    \\
    \hline
  \end{tabular}
  \begin{tabular}{|c||c|c|c|c|c|}
    \hline%\rule{0pt}{2.5ex}
    \multirow{2}{*}{$\tau$}  
    & \multicolumn{5}{|c|}{\rule{0pt}{2.5ex}$L^2$-error\qquad $\eps=10^{-20}$} \\
    \cline{2-6} 
    &\rule{0pt}{2.5ex}
    P& $E_{AP}$& $E_{APS}$& $RK_{AP}$& $RK_{APS}$\\
    \hline
    \hline\rule{0pt}{2.5ex}
    0.1 &
    $2.28\times 10^{-1} $&
    $ 1.14\times 10^{-3} $&
    $ 1.14\times 10^{-3} $&
    $ 7.4\times 10^{-5} $&
    $ 7.4\times 10^{-5} $
    \\
    \hline\rule{0pt}{2.5ex}
    0.05 &
    $2.53\times 10^{-1} $&
    $6.2\times 10^{-4} $&
    $6.2\times 10^{-4} $&
    $2.10\times 10^{-5} $&
    $2.10\times 10^{-5} $
    \\
    \hline\rule{0pt}{2.5ex}
    0.025 &
    $2.53\times 10^{-1}$&
    $3.07\times 10^{-4} $&
    $3.07\times 10^{-4} $&
    $5.3\times 10^{-6} $&
    $5.3\times 10^{-6} $
    \\
    \hline\rule{0pt}{2.5ex}
    0.0125 &
    $2.50\times 10^{-1} $&
    $1.51\times 10^{-4} $&
    $1.51\times 10^{-4} $&
    $1.33\times 10^{-6} $&
    $1.33\times 10^{-6} $
    \\
    \hline\rule{0pt}{2.5ex}
    0.00625 &
    $2.51\times 10^{-1} $&
    $7.6\times 10^{-5} $&
    $7.6\times 10^{-5} $&
    $3.44\times 10^{-7} $&
    $3.44\times 10^{-7} $
    \\
    \hline\rule{0pt}{2.5ex}
    0.003125 &
    $2.53\times 10^{-1} $&
    $3.81\times 10^{-5} $&
    $3.81\times 10^{-5} $&
    $1.18\times 10^{-7} $&
    $1.18\times 10^{-7} $
    \\
    \hline\rule{0pt}{2.5ex}
    0.0015625 &
    $2.53\times 10^{-1} $&
    $1.90\times 10^{-5} $&
    $1.90\times 10^{-5} $&
    $8.3\times 10^{-8} $&
    $8.3\times 10^{-8} $
    \\
    \hline
  \end{tabular}
  \caption{The absolute error of $u$ in the $L^{2}$-norm 
    for different time step using the singular
    perturbation scheme (P) and two proposed AP-schemes for mesh size
    $200\times 200$ at time $t=0.1$.
  }
  \label{tab:conv_t}
\end{table}

\def\xxxa{0.45\textwidth}
\begin{figure}[!ht]
  \centering
  \subfigure[$\eps=1$]
  {\includegraphics[angle=0,width=\xxxa]{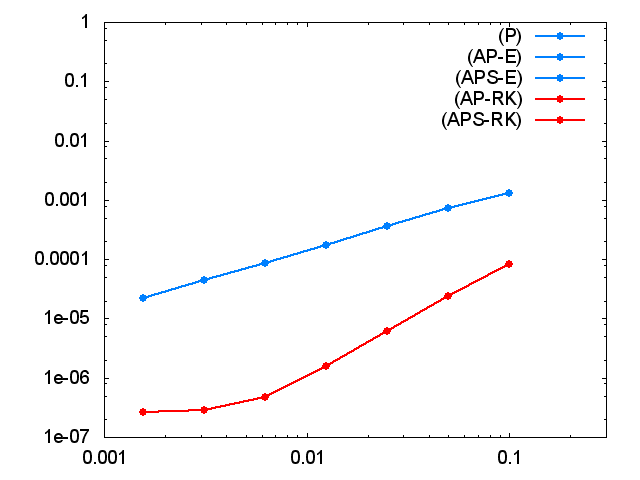}}
  \subfigure[$\eps=10^{-20}$]
  {\includegraphics[angle=0,width=\xxxa]{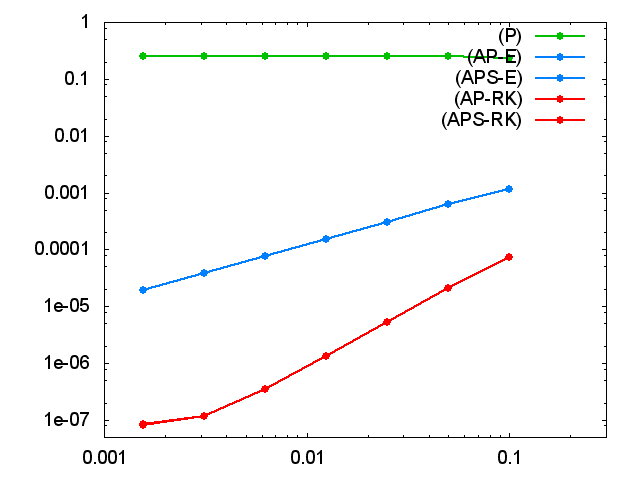}}

  \caption{Relative $L^{2}$-errors between the exact solution
    $u^{\varepsilon }$ and the computed solution with the standard
    scheme (P), the Euler-AP method ($E_{AP}$) the stabilized
    Euler-APS method ($E_{APS}$), the DIRK-AP scheme ($RK_{AP}$) and
    the stabilized DIRK-APS scheme ($RK_{APS}$) as a function of
    $\tau$, and for $\eps=1$ resp. $\eps=10^{-20}$ and a mesh with
    $200\times 200$ points.}
  \label{fig:conv_t2}
\end{figure}

%% To conclude, one can remark that the asymptotic-preserving schemes,
%% ($E_{AP}$) and ($RK_{AP}$), are uniformly accurate with respect to the
%% perturbation parameter $\eps$. This essential feature can be very
%% useful in situations where the anisotropy  is variable in space, {\it
  %% i.e.} the parameter $\eps(x)$ is $x$-dependent. No mesh-adaptation
%% is any more needed in these cases, a simple Cartesian grid enables
%% accurate results, with no regard to the $\eps$-values.

Numerical test have confirmed that all asymptotic preserving schemes
are accurate independently of $\eps$. The penalty stabilization
procedure can introduce a very small amount of additional error in
some configurations, but the order of convergence is conserved. In the
following experiments, the stabilized schemes are tested on the case
of $b$ containing closed field lines.

\subsection{Temperature balance in the presence of magnetic islands}

In this section we perform a numerical experiment related to the
tokamak plasma. This numerical test case fully demonstrates the
novelty of the proposed stabilized AP scheme as it can applied in more
general settings as the previous ones \cite{DLNN,LNN}. The results of
this section show that it is capable of simulating the heat transfer
in the presence of the closed field lines in the computational domain.

We consider a square computational domain $\Omega =[-0.5 , 0.5] \times
[-0.5 , 0.5] $ and a field $b$ with a perturbation consisting of a
region with closed lines. This so called magnetic island is initially
localized in the center of the domain and moving with the velocity
$\omega $. The field is given by
\begin{gather}
  b = \frac{B}{|B|}\, , \quad
  B =
  \left(
  \begin{array}{c}
    -A 2\pi \sin(2\pi (y-\omega t)) \nonumber\\%\cos (\pi x) \nonumber\\
    \pi \sin (\pi x) %- A\pi \cos(2\pi (y-\omega t))\sin (\pi x) 
  \end{array}
  \right)
%  \label{eq:Jfdb}
\nonumber 
,
\end{gather}
where $A$ is a perturbation parameter related to the island's width $w
= 4A^{1/2}/\pi $.  This is the largest distance between the two
branches of the separatrix, the line that divides the domain into
regions of open and closed field lines. The two branches meet at the
X-point, the saddle point of the vector potential.  The island center,
an extremum of the vector potential, is referred to as the O-point.
If $A=0$ the obtained field is aligned with the $Y$ axis and points
upwards (downwards) for $x>0$ ($x<0$). For $A>0$ the magnetic island
consisting of closed field lines appears in the region around $x=0$.
 
Referring to the theory of magnetic islands, this magnetic field
geometry approximates a saturated tearing mode in the so called
constant-$\psi$ regime (whence the parameter $A$ is a constant). The
frequency describes the rotation of an island as observed in
experiments and numerical simulations.
 
This frequency is typically of the order of the diamagnetic frequency,
which is smaller than the Alfven frequency (the propagation rate of an
Alfven wave), but larger than the transport rate across a typical
island.  It is thus interesting to study heat transport in an island
that rotates sufficiently fast.

For a static island or one rotating sufficiently slowly, we expect the
fast transport along the field lines to flatten the temperature
profile in the island region.

We choose to impose periodic boundary conditions on $\{ (x,y) \in
\partial\Omega \ | y = -0.5 \} \cup \{ (x,y) \in \partial\Omega \ | y
= 0.5 \}$ so the domain is topologically equivalent to the surface of
a torus.  We supply the computational domain with two sets of
conditions on remaining boundaries:
\begin{enumerate}

\item Dirichlet boundary conditions $\Gamma _D = \{ (x,y) \in \partial\Omega \ |
  x = -0.5 \} \cup \{ (x,y) \in \partial\Omega \ |
  x = 0.5 \}$ with
  \begin{gather}
    g_D(t,\cdot) = \left\{
    \begin{array}{ll}
      1 \quad &\text{ for } x=-0.5 \nonumber\\
      0 \quad &\text{ for } x=0.5 
    \end{array}
    \right.
    \label{eq:bcD},
  \end{gather}
  where temperature is exchanged with the exterior only by a $\Gamma_D$ boundary.

\item Neumann and Dirichlet boundary conditions: $\Gamma _N = \{ (x,y) \in \partial\Omega \ |
  x = -0.5 \}$ and $\Gamma _D = \{ (x,y) \in
  \partial\Omega \ | x = 0.5 \}$ with 
  \begin{gather}
    \qquad
    g_N(t,\cdot) = 1,
    \qquad\qquad
    g_D(t,\cdot) = 0
    \label{eq:bcN},
  \end{gather}
  which corresponds to the constant heating of the left side of the
  computational domain. 

\end{enumerate}
In the case of the boundary conditions of the first type we expect
that the presence of the island should increase the gradient of the
temperature outside the island region and keep the temperature
constant inside the island in such way that the total energy of the
system remains unchanged. If the boundary conditions are of the second
type, {\it i.e.} in the heating case, the gradient of the temperature
in the non-island region should remain constant leading thus to a loss
of the total energy of the system.

We perform simulations with fixed $\eps = 10^{-10}$ and $A=0.01$,
giving an island of a width $w=0.4/\pi \approx 13 \% $. The magnetic field lines are
shown on Figure \ref{fig:mi}. The initial conditions for both boundary
types are the same and correspond to the stationary solution with no
island present. That is to say, $u^0(x,y) = -x + \frac{1}{2} $. We
perform our simulations on a fixed grid of $200\times 200$ points for
100 time steps $\tau =2.5\times 10^{-3}$ until the final time $0.25$
is reached. We compare results for a stationary island ($\omega
=0$) and a moving island ($\omega = 10$) with no island case. We are
interested in the temperature profile along the $X$ axis as well as in
the total energy of the system, {\it i.e.} integral of the temperature
in the computational domain.

\def\xxxa{0.45\textwidth}
\begin{figure}[!t]
  \centering
  {\includegraphics[angle=0,width=\xxxa]{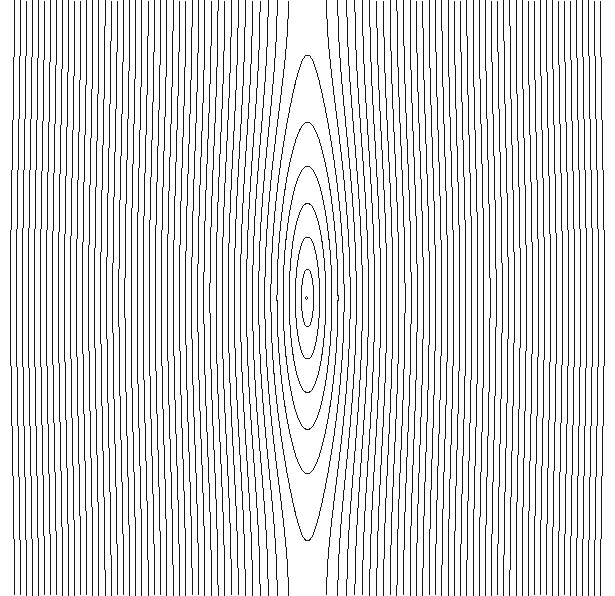}}
  \caption{Magnetic island for $A=0.01$}
  \label{fig:mi}
\end{figure}

\subsubsection{Dirichlet boundary condition on the left edge
  (\ref{eq:bcD})}

\def\xxxa{0.35\textwidth}
\begin{figure}[!ht]
  \begin{center}
    \begin{tabular}{cc}
      {\includegraphics[angle=0,width=\xxxa]{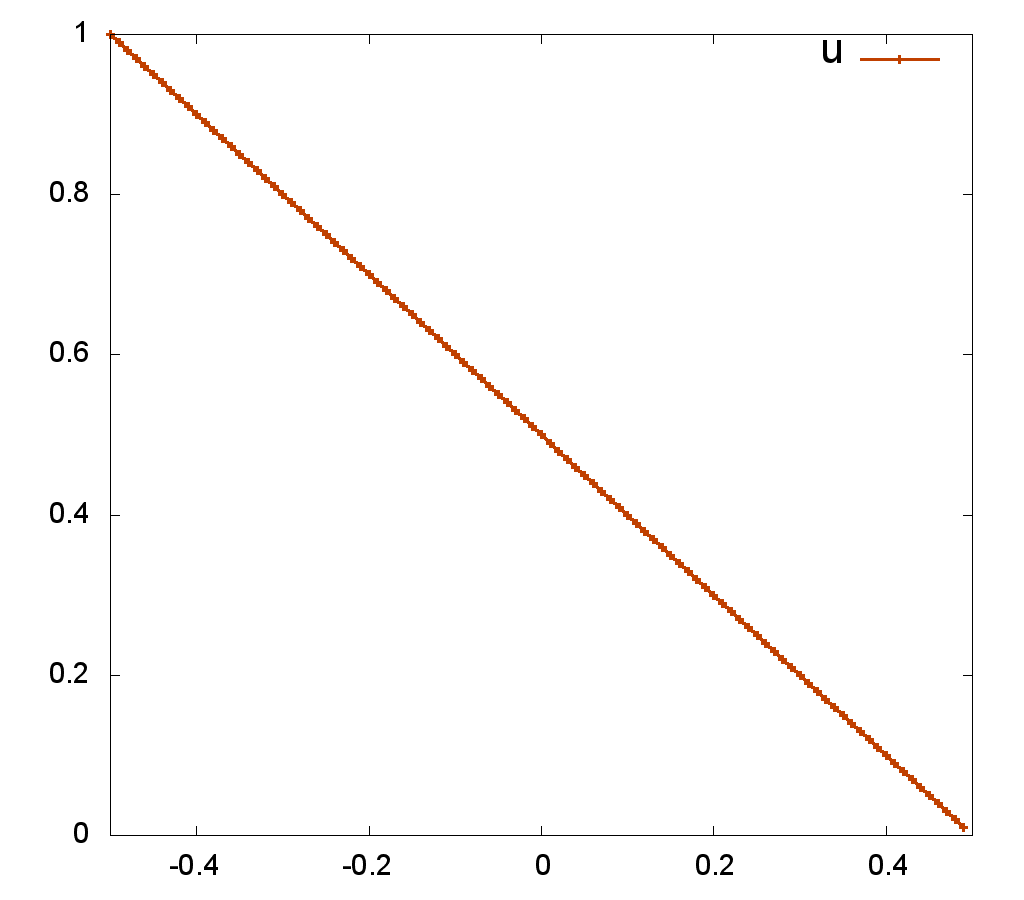}} &
      {\includegraphics[angle=0,width=\xxxa]{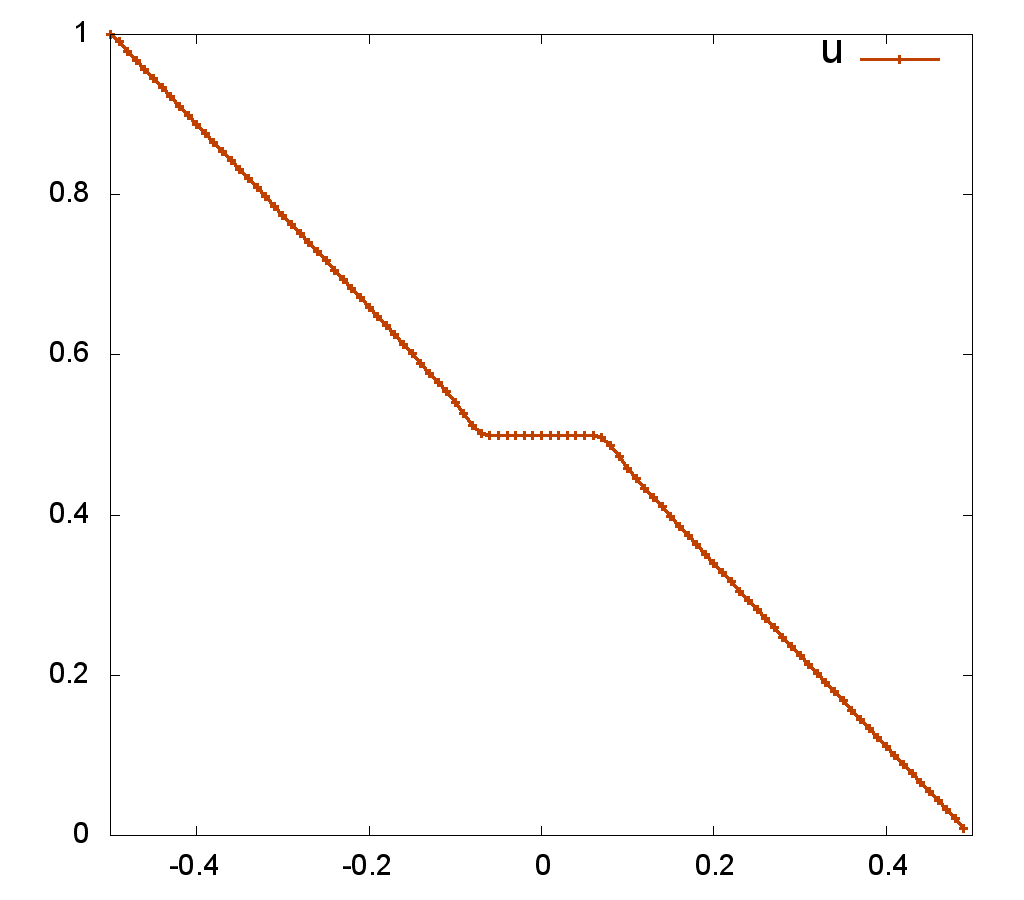}} \\
      {\includegraphics[angle=0,width=\xxxa]{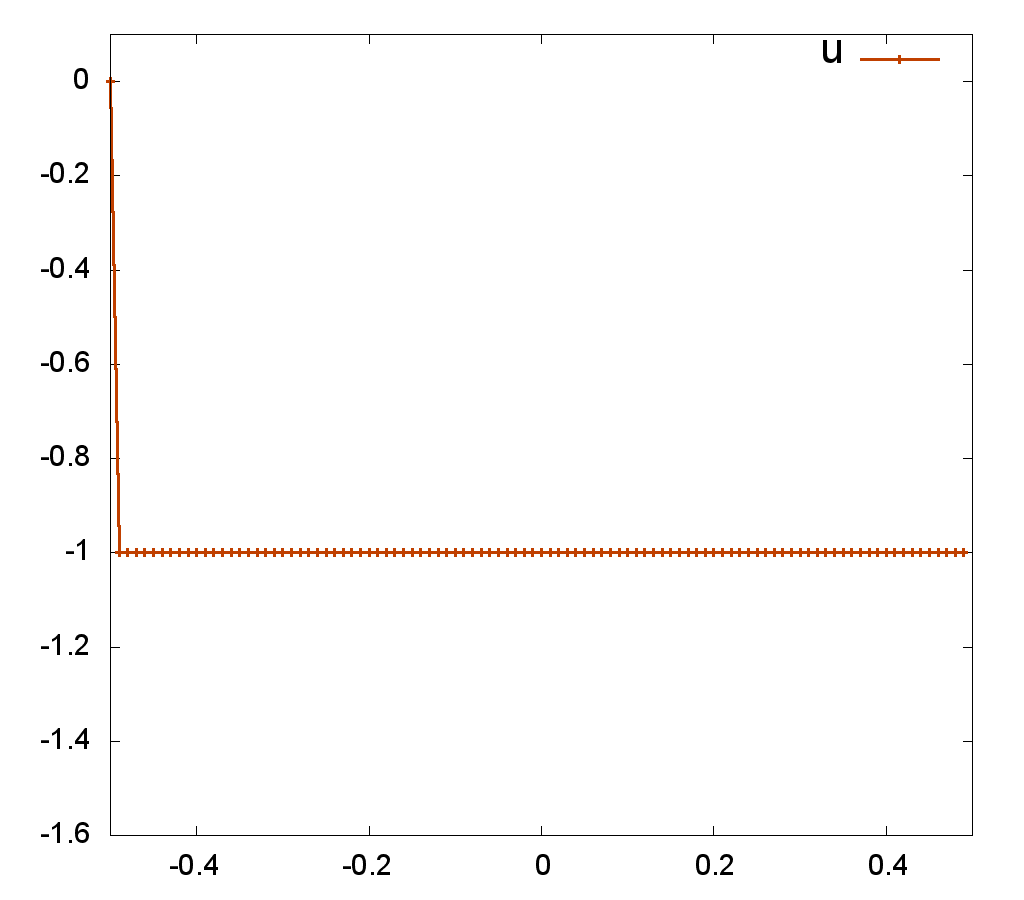}} &
      {\includegraphics[angle=0,width=\xxxa]{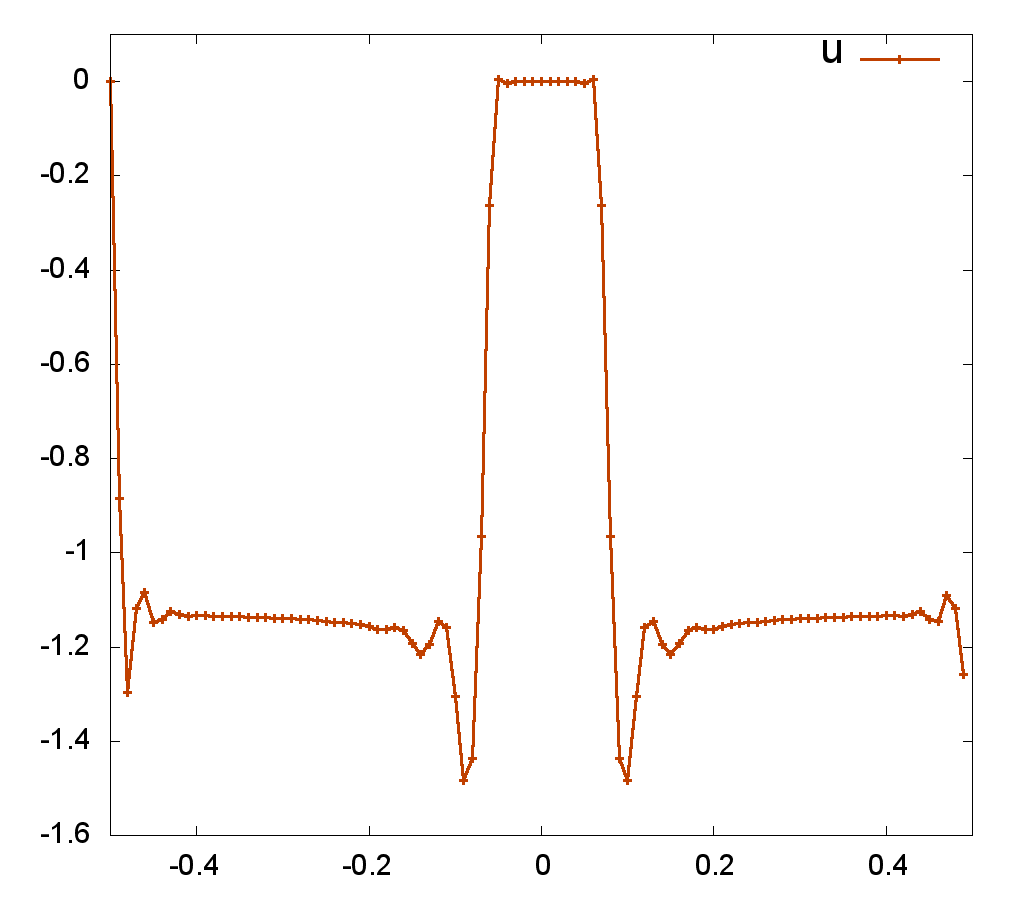}}\\
    \end{tabular}
  \end{center}
  \caption{Temperature profiles (top row) and temperature gradients
    (bottom row) along the $X$ axis for non perturbed field ($A=0$) on
    the left and a stationary island ($A=0.01$) present in the center
    of the domain on the right for the Dirichlet BC.}
  \label{fig:mi_d1}
\end{figure}

Numerical results confirm our expectations. The total energy remains
constant in the system. Integral of the temperature in the
computational domain equals to $1/2$ in all three cases: for a
stationary island, moving island and non perturbed system. The
temperature is constant in the island region leading to a stronger
gradient outside the perturbation. Temperature profiles along the $X$
axis are shown on Figures \ref{fig:mi_d1} and \ref{fig:mi_d2}. In the
case of a moving island the width of a constant temperature region in
the temperature profile along the $X$ axis is oscillating in time as
the island moves in the domain. It is interesting to note that even at
the time when there are no closed field lines across the $X$ axis the
flat region can be observed. In fact, the $x$ component of temperature
gradient is negative but close to zero in this region. Remark also
that the temperature gradient increases, approaching zero, near the
island at its largest width, i.e. across the O-point.  One refers to
this effect as "profile flattening".

\def\xxxa{0.35\textwidth}
\def\xxxb{0.3\textwidth}
\begin{figure}[!ht]
  \begin{center}
    \begin{tabular}{cc}
      {\includegraphics[angle=0,width=\xxxa]{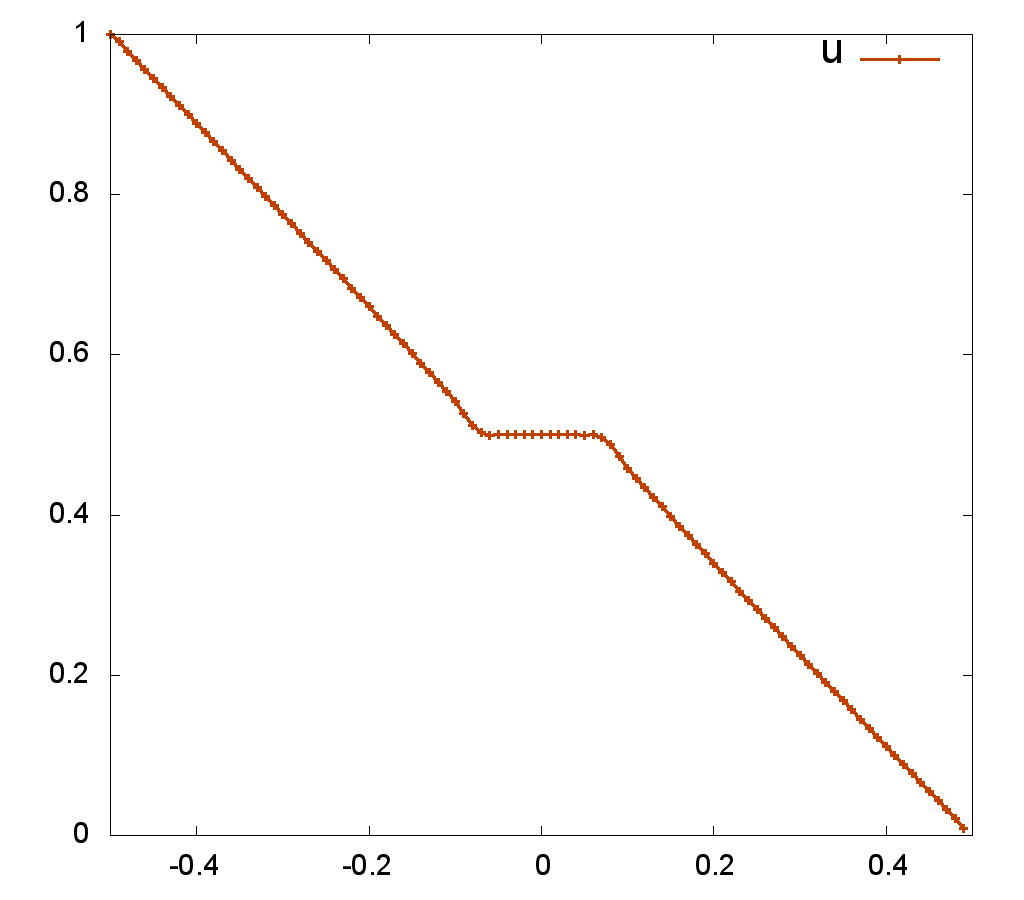}} &
      {\includegraphics[angle=0,width=\xxxa]{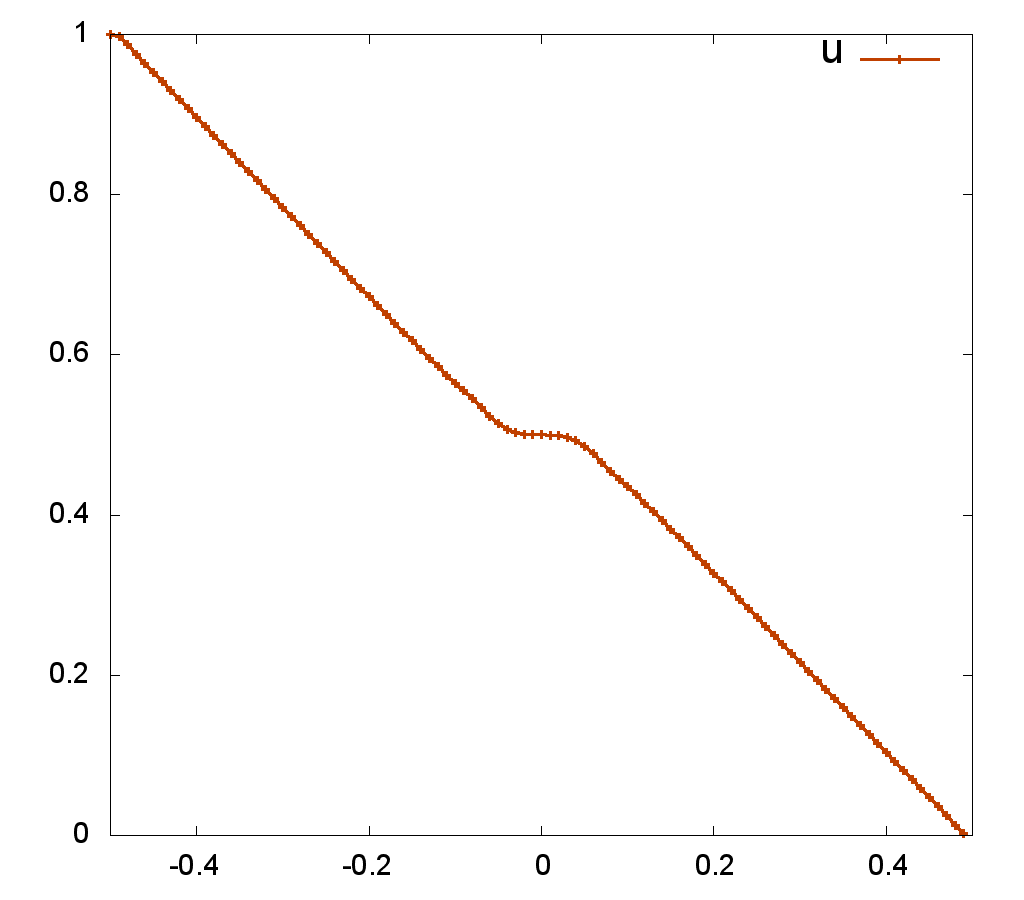}} \\
      %% \multicolumn{2}{c}{} \\
      {\includegraphics[angle=0,width=\xxxa]{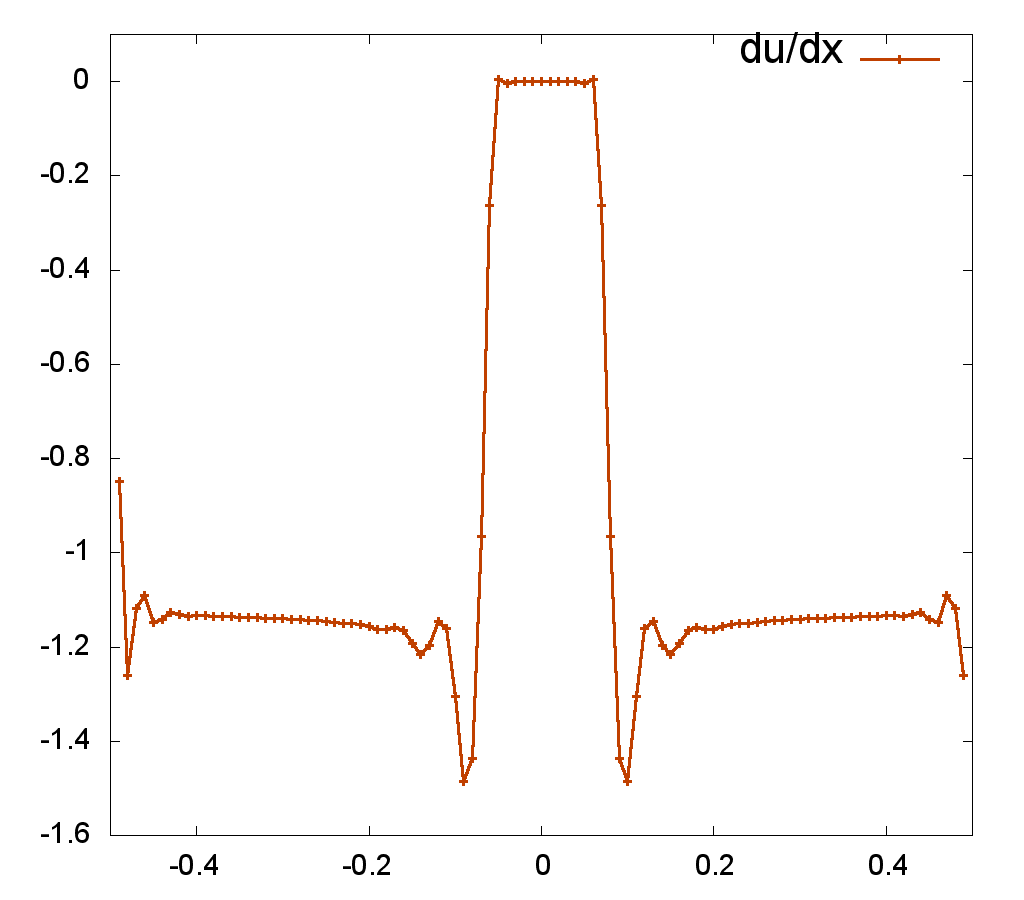}} &
      {\includegraphics[angle=0,width=\xxxa]{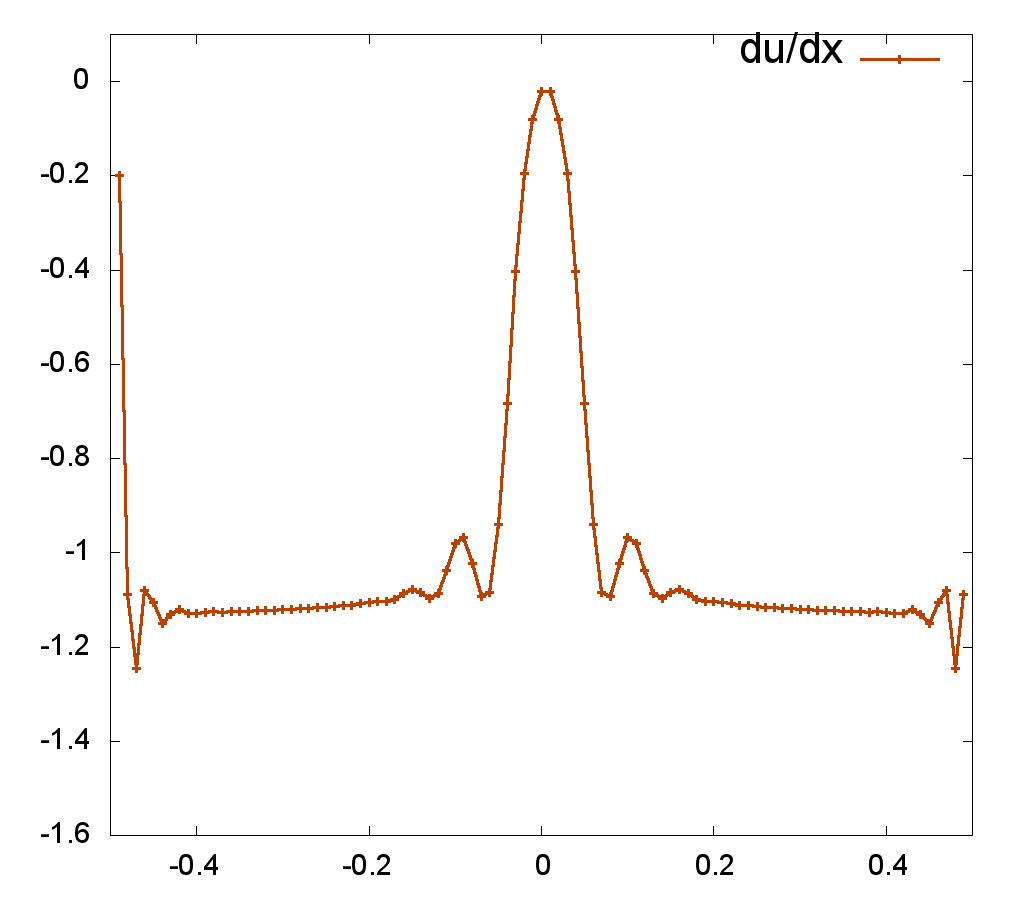}} \\
      {\includegraphics[angle=0,width=\xxxb]{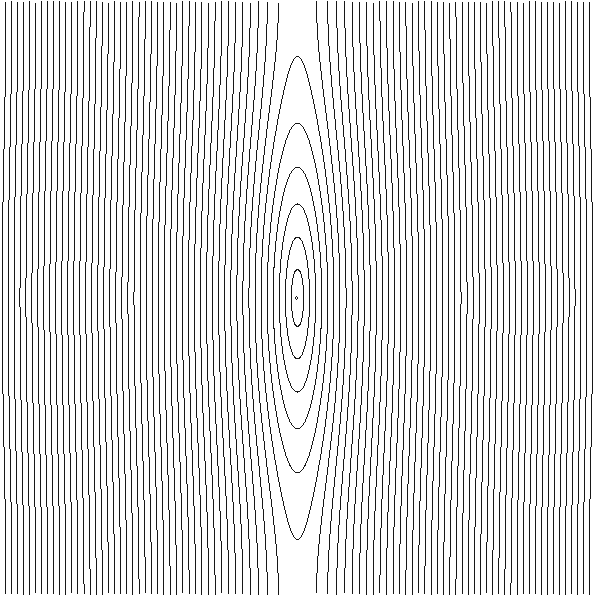}} &
      {\includegraphics[angle=0,width=\xxxb]{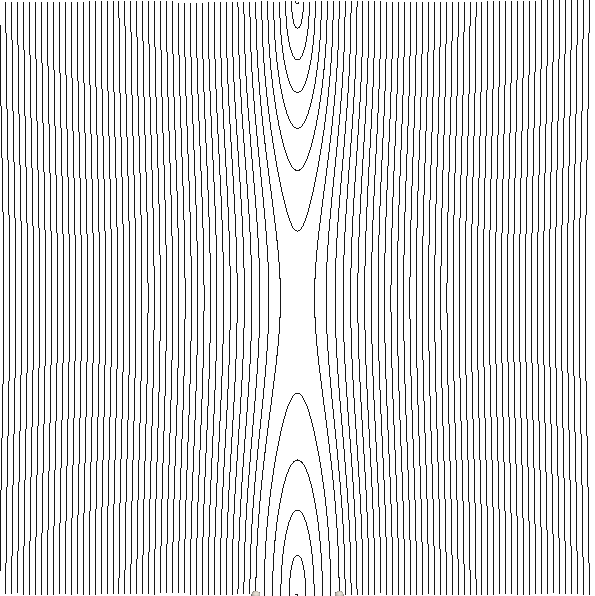}} \\
      $t=0.1$& $t=0.15$
    \end{tabular}
  \end{center}
  \caption{Temperature profiles along the $X$ axis for a moving island
    ($A=0.01$, $\omega =10$) in the first row, $x$ component of
    temperature gradient in the second row and a corresponding
    anisotropy field in the last row for different time steps for the
    Dirichlet BC.
  }
  \label{fig:mi_d2}
\end{figure}

\subsubsection{Neumann boundary condition on the left edge
  (\ref{eq:bcN})}

In the case of Neumann boundary condition imposed on the left boundary
of the domain the results again are consistent with our
expectation. The presence of the island reduces the total
energy. Integral of the temperature drops from $0.5$ to $0.44$ and the
maximal temperature from $1$ to $0.89$ for both stationary and moving
island. Temperature profiles along the $X$
axis are shown on Figures \ref{fig:mi_n1} and \ref{fig:mi_n2}. The
results are very similar to those in the
previous test case. The difference is in the maximal temperature and
in the temperature gradient in the $X$ direction, which is now close
to one in the regions far from the island. This is of course what is
expected and consistent with the heating imposed on the left boundary
by the means of Neumann boundary condition.

\def\xxxa{0.35\textwidth}
\begin{figure}[!ht]
  \begin{center}
    \begin{tabular}{cc}
      {\includegraphics[angle=0,width=\xxxa]{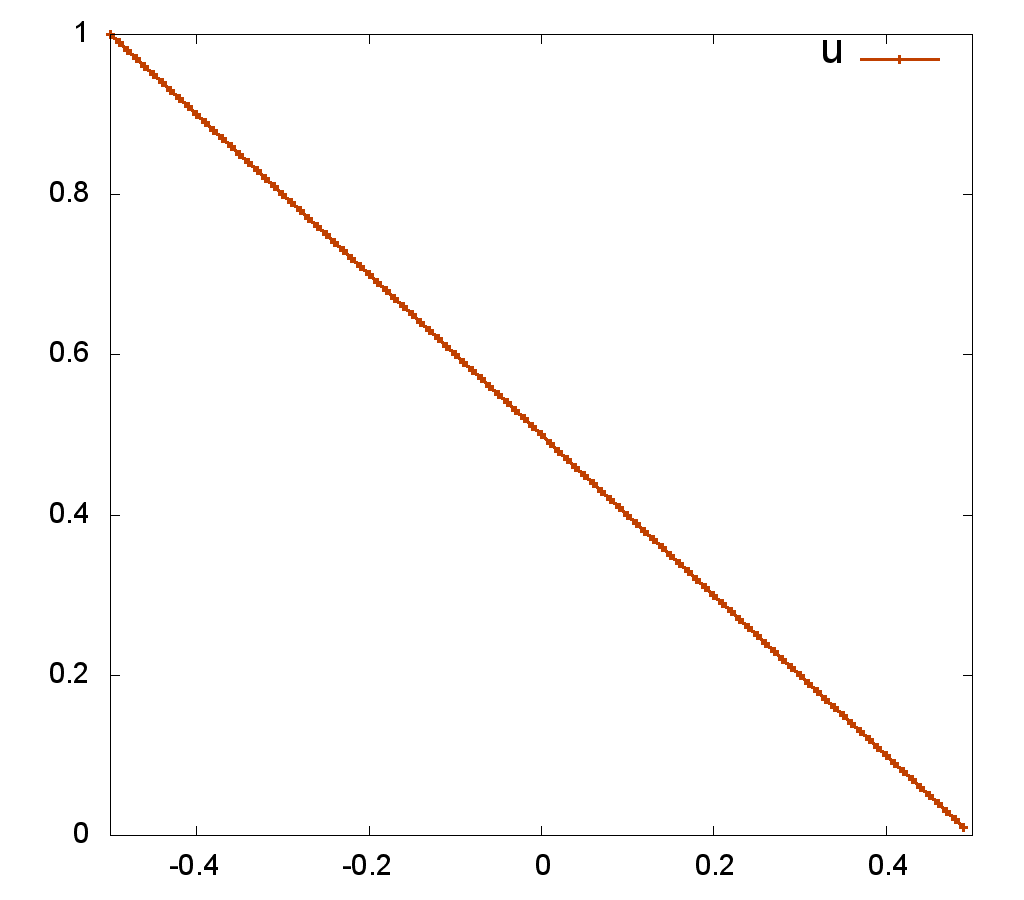}} &
      {\includegraphics[angle=0,width=\xxxa]{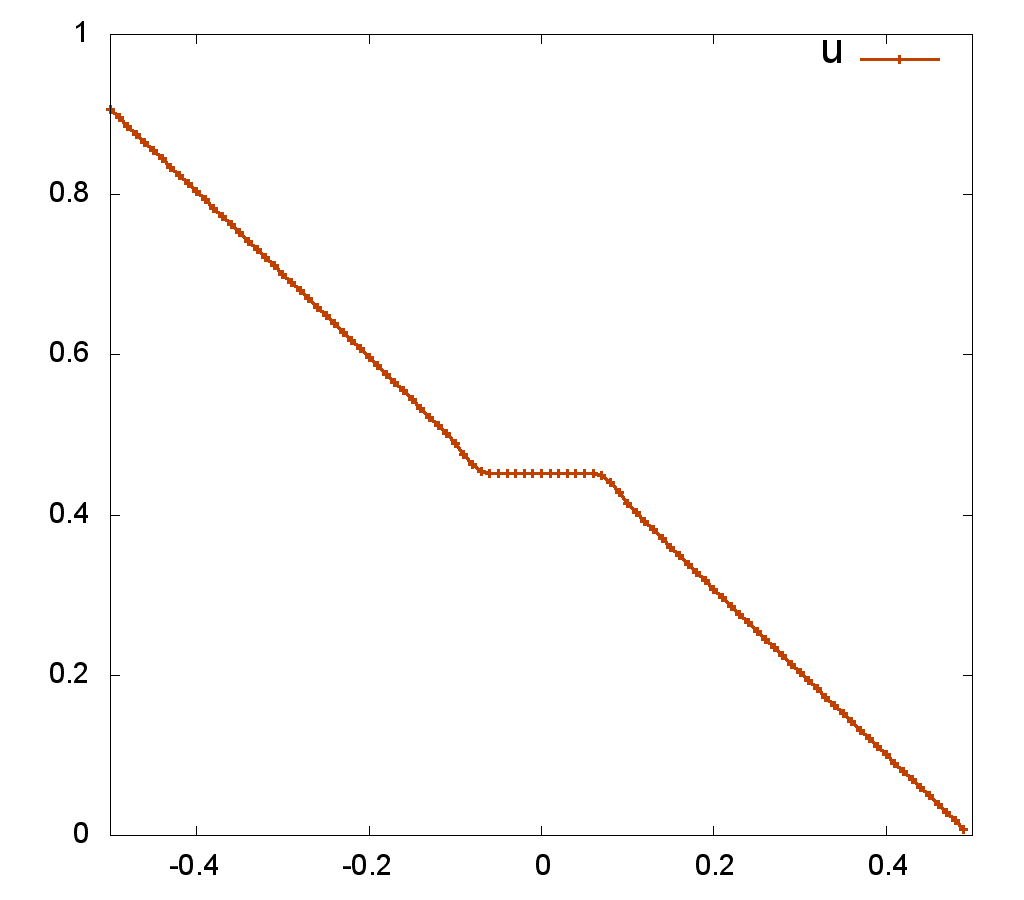}} \\
      {\includegraphics[angle=0,width=\xxxa]{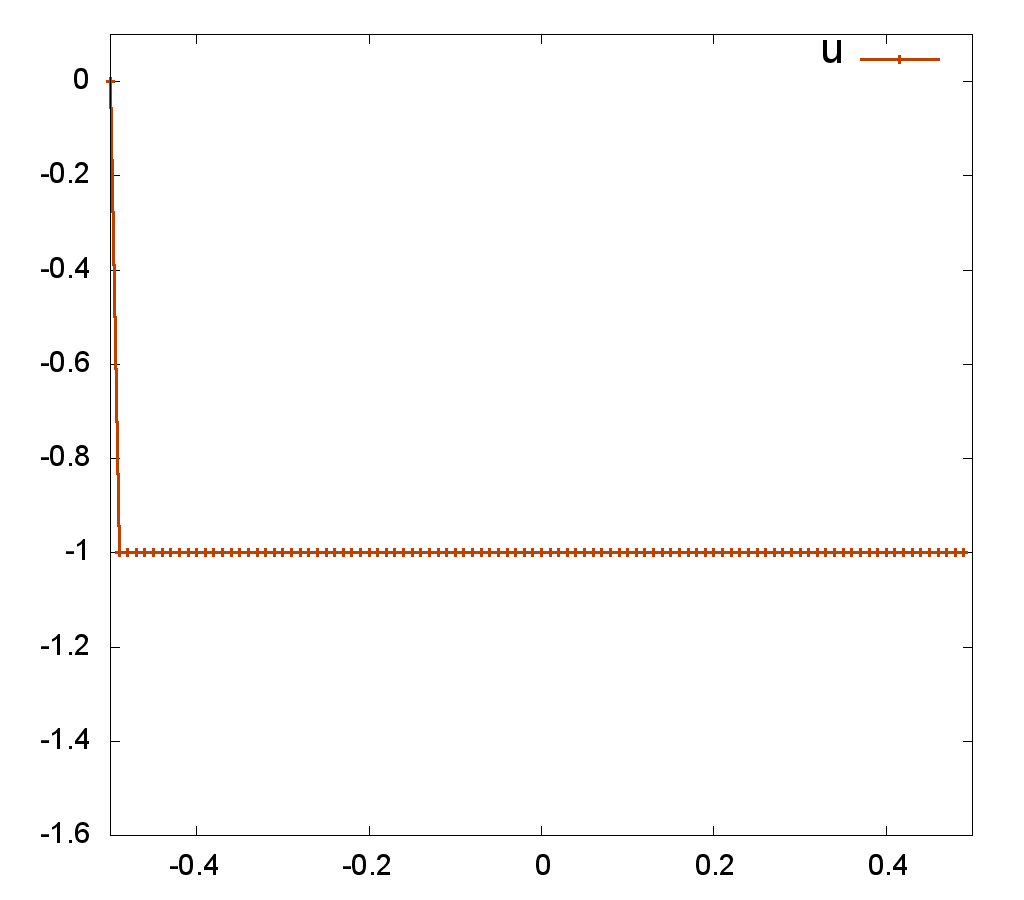}} &
      {\includegraphics[angle=0,width=\xxxa]{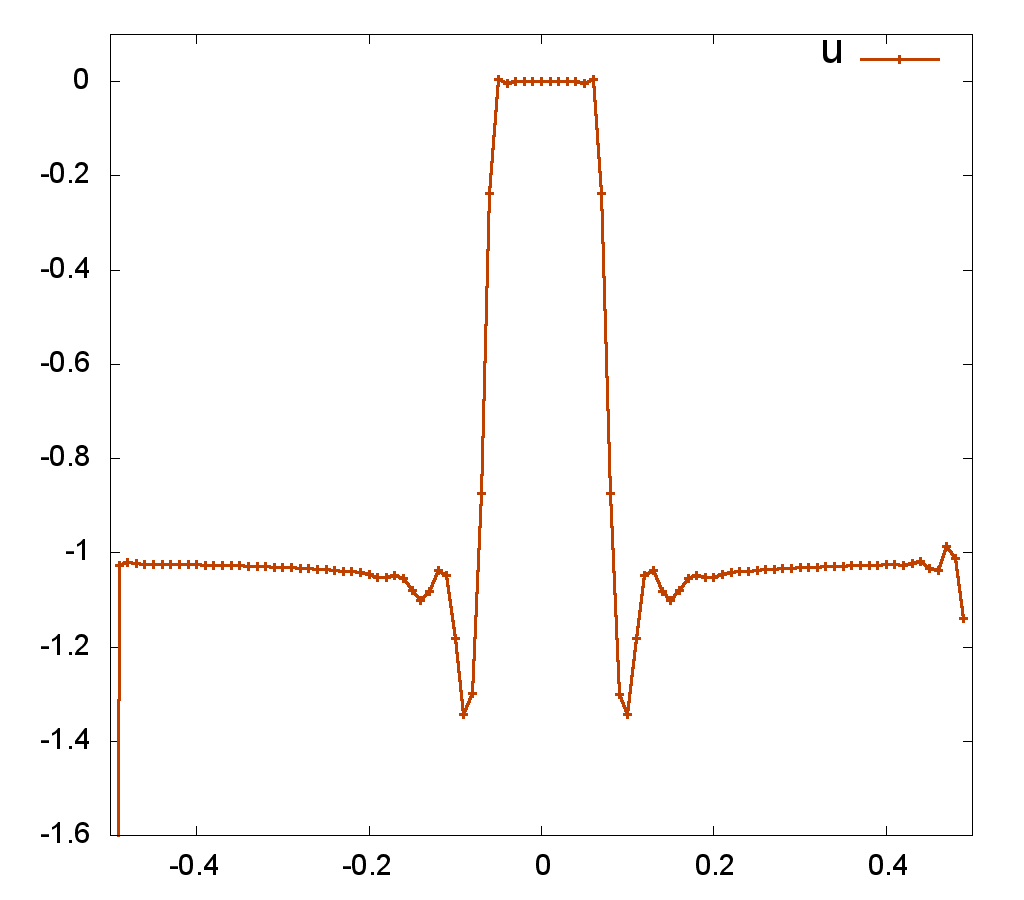}}\\
    \end{tabular}
  \end{center}
  \caption{Temperature profiles (top row) and temperature gradients
    (bottom row) along the $X$ axis for non perturbed field ($A=0$) on
    the left and a stationary island ($A=0.01$) present in the center
    of the domain on the right for the Neumann BC.}
   \label{fig:mi_n1}
\end{figure}
\def\xxxa{0.35\textwidth}
\begin{figure}[!ht]
  \begin{center}
    \begin{tabular}{cc}
      {\includegraphics[angle=0,width=\xxxa]{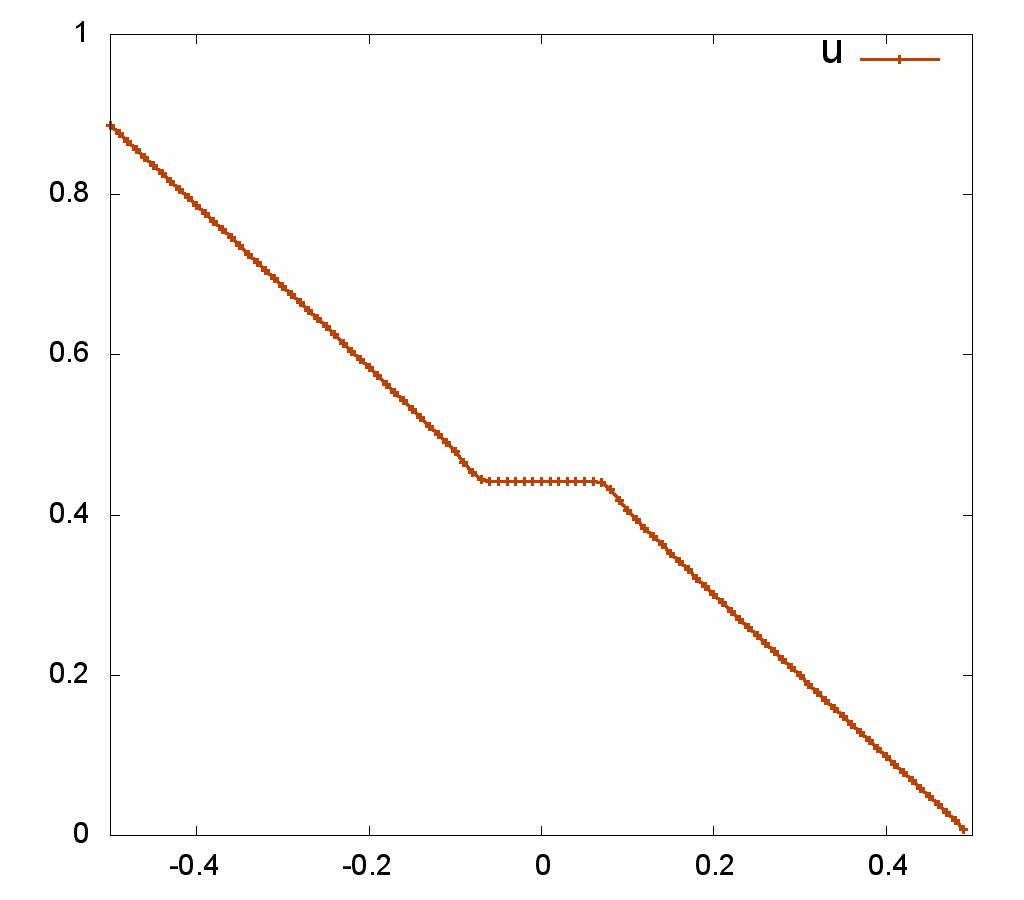}} &
      {\includegraphics[angle=0,width=\xxxa]{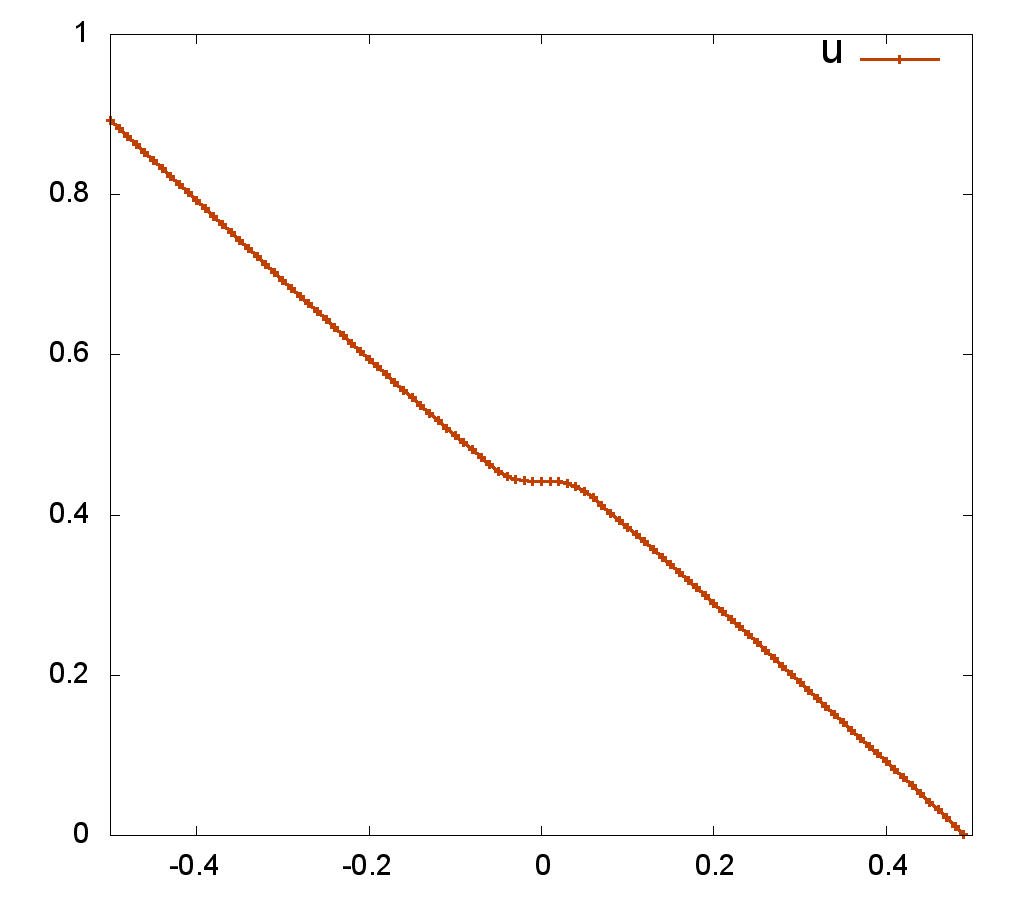}} \\
      %% \multicolumn{2}{c}{} \\
      {\includegraphics[angle=0,width=\xxxa]{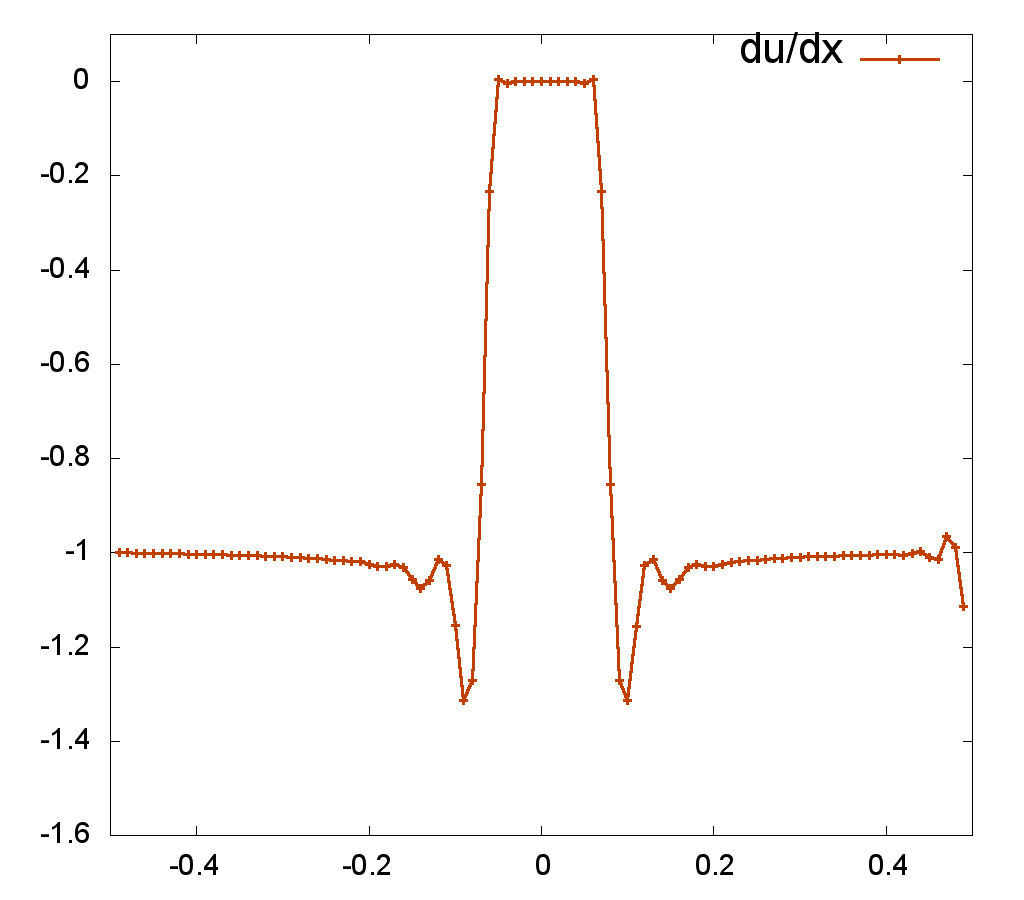}} &
      {\includegraphics[angle=0,width=\xxxa]{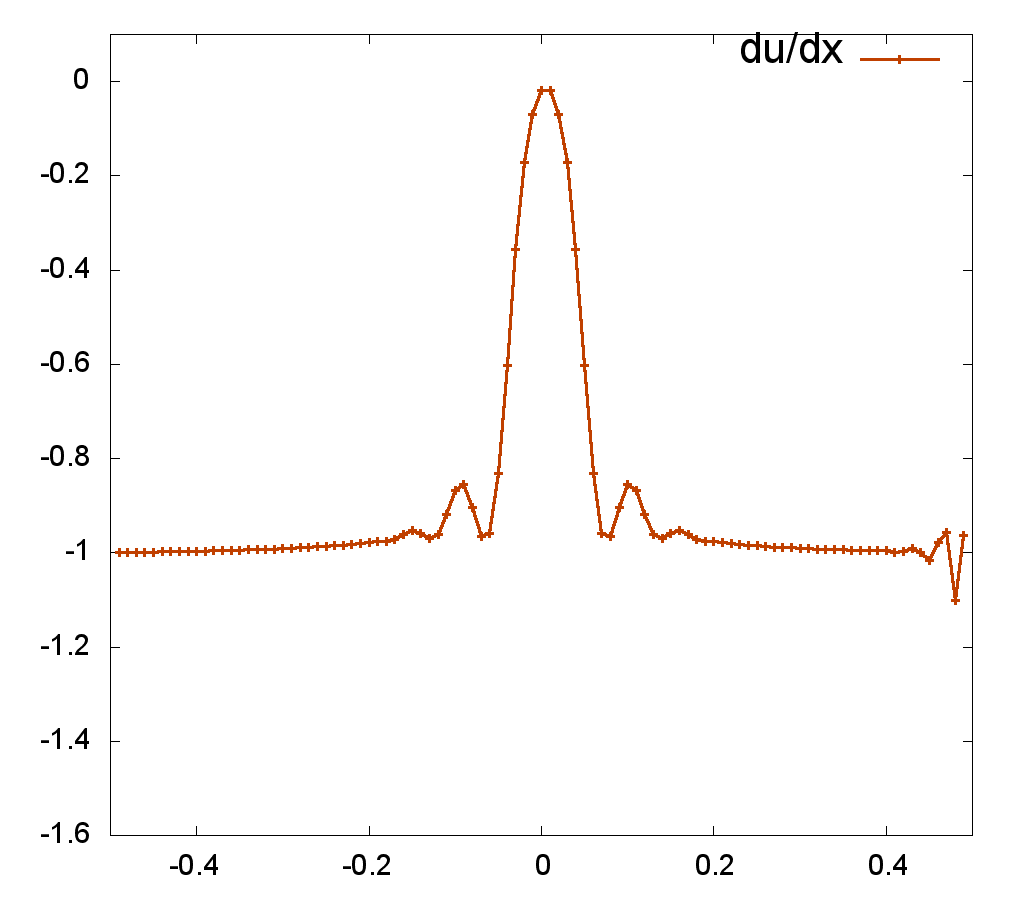}} \\
      $t=0.1$& $t=0.15$
    \end{tabular}
  \end{center}
  \caption{Temperature profiles along the $X$ axis for a moving island
    ($A=0.01$, $\omega =10$) in the first row and the $x$ component of
    temperature gradient in the second row for different time steps
    for the Neumann BC.
  }
  \label{fig:mi_n2}
\end{figure}

\subsubsection{Fast rotating magnetic island with Dirichlet boundary condition on the left edge
  (\ref{eq:bcD})}

In the last experiment we study the effect of the islands rotation
speed on the temperature profile. We expect to see a different
temperature profile when the rotation is faster than the transport
rate in the parallel direction. To achieve this numerically we
decrease the island width such that $A=0.000625$, augment the mesh size
to $500 \times 500$, decrease the size of the computational domain
$\Omega =[-0.125 , 0.125] \times [-0.5 , 0.5] $ and put
$\eps=10^{-3}$. We also increase the time resolution and put
$\tau=0.25 \times 10^{-6}$ and vary the rotation speed from $10^{3}$
to $10^{6}$. Then we compare temperature profiles after 10000 time
steps.

\def\xxxa{0.45\textwidth}
\begin{figure}[!ht]
  \centering
  \subfigure[]
  {\includegraphics[angle=0,width=\xxxa]{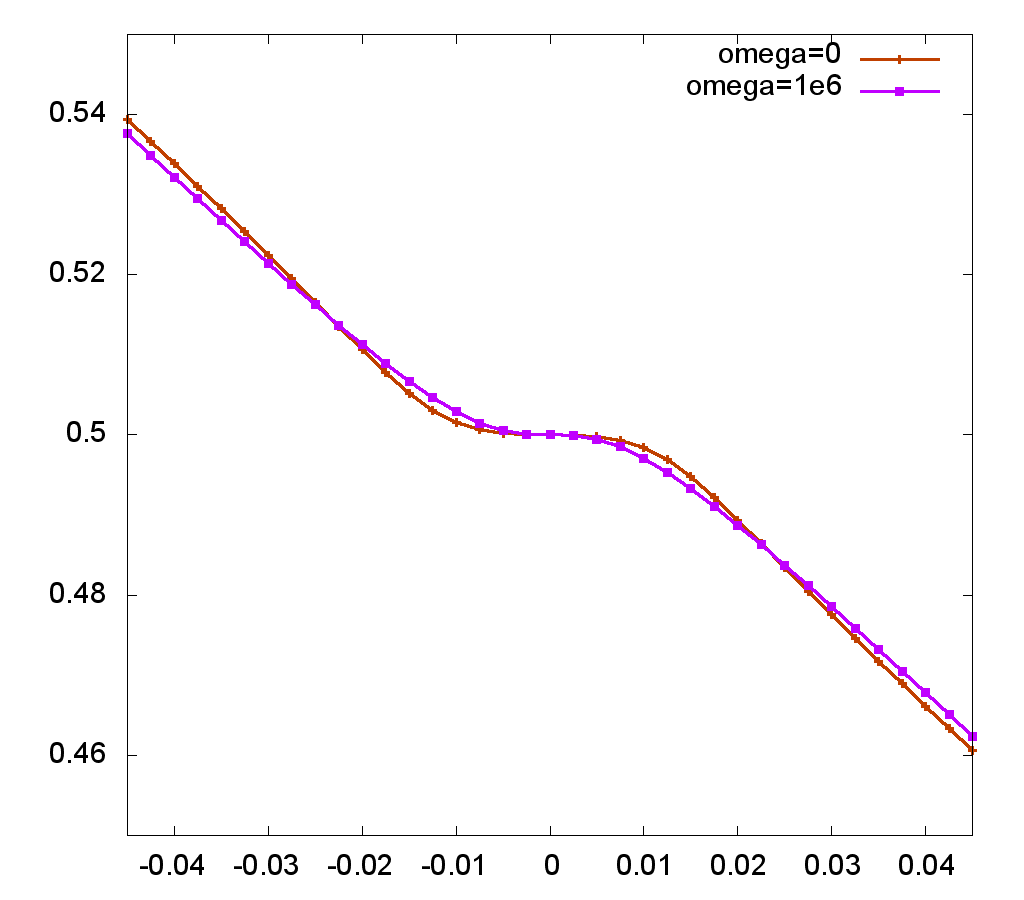}}
  \subfigure[]
  {\includegraphics[angle=0,width=\xxxa]{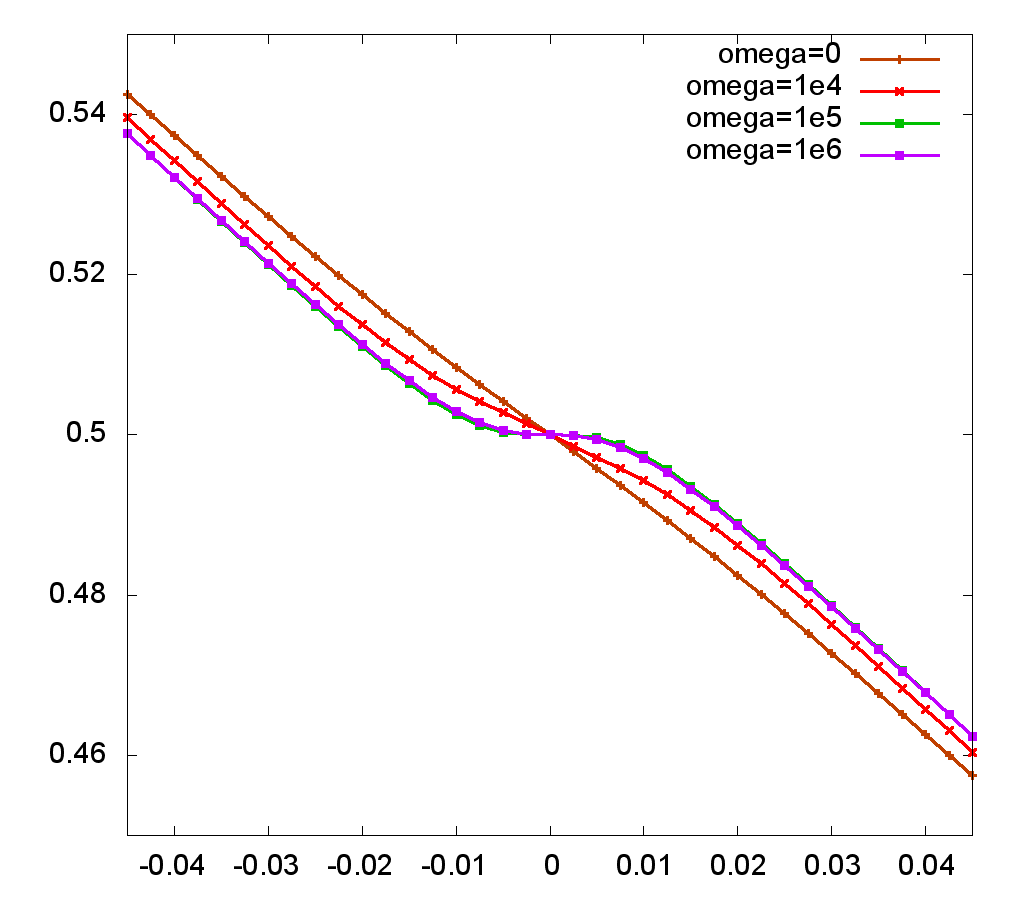}}
  \caption{Comparison between temperature profiles
    along the $X$ axis for the stationary and rotating islands.
    The profile across the island center ($Y=Y_C$) on the left 
    and away from the center ($Y=Y_C \pm 0.5$) on the right.}
  \label{fig:mi_zoom}
\end{figure}

For the smallest rotation velocity ($\omega =10^3$) no deviation from
the stationary case is observed. For other velocities the temperature
profile across the islands center is only slightly affected by the
rotation velocity. The most interesting effect of the velocity is
visible on the profile taken away from the island. In the stationary
case the profile is a straight line with a constant gradient. If the
rotation velocity is big enough, the profile begins to slightly
flatten near the center ($X=0$) for $\omega =10^4$ and finally looks
exactly the same as for the islands center ($\omega=10^5$ and
$\omega=10^6$). In fact, for a sufficiently big rotation speed the
temperature profile becomes homogeneous, {\it i.e.} independent of
$Y$. A zoom of temperature profiles at the islands center ($Y=Y_C$)
and at the most distant from the islands center $Y$ ($Y=Y_C +0.5$ if
$Y_C <0$ and $Y=Y_C-0.5$ otherwise) is presented on Figure
\ref{fig:mi_zoom}.

%%%%%%%%%%%%%%%%%%%%%%%%%%%%
\section{Conclusion}\label{SEC6}
%%%%%%%%%%%%%%%%%%%%%%%%%%%%

The here presented Asymptotic-Preserving scheme proves to be an
efficient, general and easy to implement numerical method for solving
strongly anisotropic parabolic problems. As numerical experiments
show, the second order two stage Runge-Kutta scheme is of particular
interest as it produces accurate results for relatively large time
steps allowing to substantially reduce the computational time.  The
herein introduced stabilization term is an important improvement to
the previously proposed Asymptotic Preserving schemes. The not so rare
in real applications anisotropy topologies can be successfully
addressed with this new method. Finally, the case of magnetic islands
studied in the last section brings new and interesting results showing
the homogenization of the temperature profile for fast rotating islands.

%%%%%%%%%%%%%%%%%%%%%%%%%%%%%%%%%%%%%%%%%%%%%%%%%%%%%%%%%%%%%%%%%%%%%%%%%%%%%%%%%%%%%%%%%%%%
\section*{Acknowledgments}

This work has been partially supported by the ANR project BOOST
(Building the future Of numerical methOdS for iTer, 2010-2014), and by
the European Communities under the contract of Association between
EURATOM and CEA. The views and opinions expressed herein do not
necessarily reflect those of the European Commission.

\newpage
\bibliographystyle{abbrv}
\bibliography{bib_aniso}

\end{document}